\documentclass{article}%
\usepackage{amssymb}
\usepackage{graphicx}
\usepackage{amsmath}
\usepackage{amsfonts}%
\newtheorem{theorem}{Theorem}

\newtheorem{definition}[theorem]{Definition}
\newtheorem{example}[theorem]{Example}

\newtheorem{lemma}[theorem]{Lemma}
\newtheorem{notation}[theorem]{Notation}

\newtheorem{proposition}[theorem]{Proposition}
\newtheorem{remark}[theorem]{Remark}

\newenvironment{proof}[1][Proof]{\textbf{#1.} }{\ \rule{0.5em}{0.5em}}
\begin{document}

\title{A Functorial Approach to the Infinitesimal Theory of Groupoids}
\author{Hirokazu Nishimura\\Institute of Mathematics, University of Tsukuba\\Tsukuba, Ibaraki, 305-8571, Japan}
\maketitle

\begin{abstract}
Lie algebroids are by no means natural as an infinitesimal counterpart of
groupoids. In this paper we propose a functorial construction called
\textit{Nishimura algebroids} for an infinitesimal counterpart of groupoids.
Nishimura algebroids, intended for differential geometry, are of the same vein
as Lawvere's functorial notion of \textit{algebraic theory} and Ehresmann's
functorial notion of theory called \textit{sketches}. We study \textit{totally
intransitive} Nishimura algebroids in detail. Finally we show that Nishimura
algebroids naturally give rise to Lie algebroids.

\end{abstract}

\section{\label{s0}Introduction}

Many mathematicians innocently believe that \textit{infinitesimalization} is
no other than \textit{linearization}. We contend that linearization is only a
tiny portion of infinitesimalization. It is true that Lie algebras are the
linearization of Lie groups, but it is by no means true that Lie algebras are
the infinitesimalization of Lie groups. The fortunate success of the theory of
Lie algebras together with their correspondence with Lie groups unfortunately
enhanced their wrong conviction and blurred what are to be really the
infinitesimalization of groups and, more generally, groupoids.

In this paper we propose, after the manners of Lawvere's functorial
construction of \textit{algebraic theory} and Ehresmann's functorial notion of
theory called \textit{sketches, }a functorial construction of
\textit{Nishimura algebroids} for the infinitesimalization of groupoids. After
giving some preliminaries and fixing notation in the coming section, we will
introduce our main notion of Nishimura algebroid in 6 steps. Then we will
study totally intransitve Nishimura algebroids, in which the main result is
that the linear part of any totally intransitve Nishimura algebroid is a Lie
algebra bundle. As our final investigation we will show that Nishimura
algebroids naturally give rise to Lie algebroids.

\section{\label{s1}Preliminaries}

\subsection{Synthetic Differential Geometry}

Our standard reference on synthetic differential geometry is Lavendhomme
\cite{l1}. In synthetic differential geometry we generally work within a good
topos. If the reader is willing to know how to get such a topos, he or she is
referred to Kock \cite{k1} or Moerdijk and Reyes \cite{mr}. We denote by
$\mathbb{R}$ the internal set of real numbers, which is endowed with a
cornucopia of nilpotent infinitesimals persuant to the general Kock-Lawvere
axiom. The internal category $\mathbf{Inf}$ of \textit{infinitesimal spaces}
comes contravariantly from the external category of Weil algebras over the set
of real numbers by taking $\mathrm{Spec}_{\mathbb{R}}$. In particular, the
infinitesimal space corresponding to the set of real numbers as a Weil algebra
is denoted by $1$.We should note that every infinitesimal space $\mathcal{D}$
has a distinguished point, namely, $0_{\mathcal{D}}$ (often written simply
$0$), and every morphism in $\mathbf{Inf}$\ preserves distinguished points. An
arbitrarily chosen microlinear space $M$ shall be fixed throughout the rest of
this paper.

\subsection{\label{p1}Groupoids}

Our standard reference on groupoids is \cite{ma}. Let $\mathcal{D}$ be an
object in $\mathbf{Inf}$. Given $m\in M$ and a groupoid $G$ over $M$ with its
object inclusion map $\mathrm{id}:M\rightarrow G$ and its source and target
projections $\alpha,\beta:G\rightarrow M$, we denote by $\mathcal{A}%
_{m}^{\mathcal{D}}G$ the totality of mappings $\gamma:\mathcal{D}\rightarrow
G$ with $\gamma(0_{\mathcal{D}})=\mathrm{id}_{m}$ and $(\alpha\circ
\gamma)(d)=m$ for any $d\in\mathcal{D}$. We denote by $\mathcal{A}%
^{\mathcal{D}}G$ the set-theoretic union of $\mathcal{A}_{m}^{\mathcal{D}}G$'s
for all $m\in M$. The canonical projection $\pi:\mathcal{A}^{\mathcal{D}%
}G\rightarrow M$ is defined as is expected. The \textit{anchor} $\mathbf{a}%
_{G}^{\mathcal{D}}:\mathcal{A}^{\mathcal{D}}G\rightarrow M^{\mathcal{D}}$ is
defined to be simply
\[
\mathbf{a}_{G}^{\mathcal{D}}(\gamma)=\beta\circ\gamma
\]
for any $\gamma\in\mathcal{A}^{\mathcal{D}}G$, where $M^{\mathcal{D}}$ is the
space of mappings of $\mathcal{D}$ into $M$. We note that if the groupoid $G$
is the pair groupoid $M\times M$, then $\mathcal{A}^{\mathcal{D}}(M\times M)$
can canonically be identified with $M^{\mathcal{D}}$. We write $\mathbf{I}G$
for the inner subgroupoid of $G$, for which the reader is referred to p.14 of
\cite{ma}.

\subsection{Simplicial Spaces}

The notion of \textit{simplicial space} was discussed by Nishimura \cite{n11}
and \cite{n13}, where simplicial spaces were called \textit{simplicial
objects} in the former paper, while they were called \textit{simplicial
infinitesimal spaces} in the latter paper. \textit{Simplicial spaces} are
spaces of the form
\[
D^{m}\{\mathcal{S}\}=\{(d_{1},...,d_{m})\in D^{m}|d_{i_{1}}...d_{i_{k}%
}=0\text{ for any }(i_{1},...,i_{k})\in\mathcal{S}\}\text{,}%
\]
where $\mathcal{S}$ is a finite set of sequences $(i_{1},...,i_{k})$ of
natural numbers with $1\leq i_{1}<...<i_{k}\leq m$. By way of example, we have
$D(2)=D^{2}\{(1,2)\}$ and $D(3)=D^{3}\{(1,2),(1,3),(2,3)\}$. Given two
simplicial spaces $D^{m}\{\mathcal{S}\}$ and $D^{n}\{\mathcal{T}\}$, we define
another simplicial space $D^{m}\{\mathcal{S}\}\oplus D^{n}\{\mathcal{T}\}$ to
be
\begin{align*}
&  D^{m}\{\mathcal{S}\}\oplus D^{n}\{\mathcal{T}\}\\
&  =\{(d_{1},...,d_{m},e_{1},...,e_{n})\in D^{m+n}|d_{i_{1}}...d_{i_{k}%
}=0\text{ for any }(i_{1},...,i_{k})\in\mathcal{S}\text{, }\\
e_{j_{1}}...e_{j_{l}}  &  =0\text{ for any }(j_{1},...,j_{l})\in
\mathcal{T}\text{, }d_{i}e_{j}=0\text{ for any }1\leq i\leq m\text{ and }1\leq
j\leq n\}
\end{align*}
We denote by $\mathbf{Simp}$ the full subcategory of $\mathbf{Inf}$ whose
objects are all simplicial spaces. Obviously the category $\mathbf{Simp}$ is
closed under direct products. The category $\mathbf{Simp}$ has finite
coproducts. In particular, it has the initial object $\mathbf{1}$, which is
also the terminal object.

\section{Nishimura Algebroids}

Let $M$ be a microlinear space. We will introduce our main notion of
\textit{Nishimura algebroid} \textit{over} $M$ step by step, so that the text
is divided into six subsections.

\subsection{\textit{Nishimura Algebroids}$_{1}$}

\begin{definition}
A \textit{Nishimura algebroid}$_{1}$ \textit{over} $M$ is simply a
contravariant functor $\mathcal{A}$ from the category $\mathbf{Simp}$ of
simplicial spaces to the category $\mathbf{MLS}_{M}$ of microlinear spaces
over $M$ mapping finite coproducts in $\mathbf{Simp}$ to finite products in
$\mathbf{MLS}_{M}$.
\end{definition}

Given a simplicial space $\mathcal{D}$ in $\mathbf{Simp}$, we will usually
write $\pi:\mathcal{A}^{\mathcal{D}}\rightarrow M$ for $\mathcal{A}%
(\mathcal{D})$. In particular, we will often write $\mathcal{A}^{n}$ in place
of $\mathcal{A}^{D^{n}}$. We will simply write $\pi$ for the projection to
$M$\ in preference to such a more detailed notation as $\pi_{\mathcal{A}%
,\mathcal{D}}$, which should not cause any possible confusion. Given $m\in M$,
we write $\mathcal{A}_{m}^{\mathcal{D}}$ for $\{x\in\mathcal{A}^{\mathcal{D}%
}\mid\pi(x)=m\}$. Given a morphism $f:\mathcal{D\rightarrow D}^{\prime}$ in
$\mathbf{Simp}$, we will usually write $\mathcal{A}^{f}:\mathcal{A}%
^{\mathcal{D}^{\prime}}\rightarrow\mathcal{A}^{\mathcal{D}}$ for
$\mathcal{A}(f)$. Given $m\in M$, there is a unique element in $\mathcal{A}%
_{m}^{\mathbf{1}}$, which we denote by $\mathbf{0}_{m}^{\mathbf{1}}$. Given an
object $\mathcal{D}$ in $\mathbf{Simp}$, we define $\mathbf{0}_{m}%
^{\mathcal{D}}\in\mathcal{A}_{m}^{\mathcal{D}}$ to be
\[
\mathbf{0}_{m}^{\mathcal{D}}=\mathcal{A}^{\mathcal{D\rightarrow}\mathbf{1}%
}(\mathbf{0}_{m}^{\mathbf{1}})
\]

\begin{example}
By assigning the spac $M^{\mathcal{D}}$ of mappings from $\mathcal{D}$\ into
$M$\ to each object $\mathcal{D}$ in $\mathbf{Simp}$ and assigning
$M^{f}:M^{\mathcal{D}^{\prime}}\rightarrow M^{\mathcal{D}}$ to each morphism
$f:\mathcal{D\rightarrow D}^{\prime}$ in $\mathbf{Simp}$, we have a
\textit{Nishimura algebroid}$_{1}$ \textit{over} $M$ to be called the standard
\textit{Nishimura algebroid}$_{1}$ \textit{over} $M$ and to be denoted by
$\mathcal{S}_{M}$ or more simply by $\mathcal{S}$.
\end{example}

\begin{example}
Let $G$ be a groupoid over $M$. By assigning $\mathcal{A}^{\mathcal{D}}G$ to
each object $\mathcal{D}$ in $\mathbf{Simp}$ and assigning $\mathcal{A}%
^{f}G:\mathcal{A}^{\mathcal{D}^{\prime}}G\rightarrow\mathcal{A}^{\mathcal{D}%
}G$ to each morphism $f:\mathcal{D\rightarrow D}^{\prime}$ in $\mathbf{Simp}$,
we have a \textit{Nishimura algebroid}$_{1}$ \textit{over} $M$ to be denoted
by $\mathcal{A}G$.
\end{example}

Each $\sigma\in\mathfrak{S}_{n}$ induces a morphism $\sigma:D^{n}\rightarrow
D^{n}$ in $\mathbf{Simp}$ such that
\[
\sigma(d_{1},...,d_{n})=(d_{\sigma(1)},...,d_{\sigma(n)})
\]
for any $(d_{1},d_{2})\in D^{2}$. Given $x\in\mathcal{A}^{n}$, we will often
write $^{\sigma}x$ for $\mathcal{A}^{\sigma}(x)$. It is easy to see that
\[
^{\tau\sigma}x=^{\tau}(^{\sigma}x)
\]
for any $x\in\mathcal{A}^{n}$ and any $\sigma,\tau\in\mathfrak{S}_{n}$.

Given $x\in\mathcal{A}^{n}$ and $a\in\mathbb{R}$, we define $a\underset
{i}{\cdot}x$ ($1\leq i\leq n$) to be
\[
a\underset{i}{\cdot}x=\mathcal{A}^{((d_{1},,,,,d_{n})\in D^{n}\mapsto
(d_{1},...,d_{i-1},ad_{i},d_{i+1},...,d_{n})\in D^{n})}(x)
\]

\subsection{\textit{Nishimura Algebroids}$_{2}$}

\begin{definition}
A Nishimura algebroid$_{1}$\ $A$\ over $M$ is called a \textit{Nishimura
algebroid}$_{2}$ \textit{over} $M$ if the application of $\mathcal{A}$ to any
quasi-colimit diagram in $\mathbf{Simp}$ results in a limit diagram.
\end{definition}

\begin{remark}
The notion of \textit{Nishimura algebroid}$_{2}$ \textit{over} $M$ can be
regarded as a partial algebrization of microlinearity.
\end{remark}

\begin{example}
The standard Nishimura algebroid$_{1}$ $\mathcal{S}_{M}$ over $M$ is a
Nishimura algebroid$_{2}$ over $M$. This follows simply from our assumption
that $M$ is a microlinear space.
\end{example}

\begin{example}
Let $G$ be a groupoid over $M$. Then the Nishimura algebroid$_{1}$
$\mathcal{A}G$ over $M$ is a Nishimura algebroid$_{2}$ over $M$. This follows
simply from our assumption that $M$ and $G$ are microlinear spaces.
\end{example}

Let $\mathcal{A}$ be a Nishimura algebroid$_{2}$ over $M$. Let $m\in M$ with
$x,y\in\mathcal{A}_{m}^{1}$. By using the quasi-colimit diagram (1) of small
objects referred to in Proposition 6 (\S 2.2) of Lavendhomme \cite{l1}, there
exists a unique $z\in\mathcal{A}^{D\oplus D}$ with $\mathcal{A}^{(d\in
D\longmapsto(d,0)\in D\oplus D)}(z)=x$ and $\mathcal{A}^{(d\in D\longmapsto
(0,d)\in D\oplus D)}(z)=y$. We define $x+y$ to be $\mathcal{A}^{(d\in
D\longmapsto(d,d)\in D\oplus D)}(z)$. Given $a\in\mathbb{R}$, we define $ax$
to be $\mathcal{A}^{(d\in D\mapsto ad\in D)}(x)\in\mathcal{A}_{m}^{1}$. With
these operations we have

\begin{theorem}
\label{tn2.1}Given a Nishimura algebroid$_{2}$ $\mathcal{A}$ over $M$,
$\mathcal{A}_{m}^{1}$ is an $\mathbb{R}$-module for any $m\in M$.
\end{theorem}

\begin{proof}
The proof is essentially a familiar proof that $\mathcal{S}_{m}^{1}$ is an
$\mathbb{R}$-module, for which the reader is referred, e.g., to Lavendhomme
\cite{l1}, \S 3.1, Proposition 1. What we should do is only to reformulate the
familiar proof genuinely in terms of diagrams. The details can safely be left
to the reader.
\end{proof}

Let $\mathcal{A}$ be a Nishimura algebroid$_{2}$ over $M$ with $m\in M$.

Let $x,y\in\mathcal{A}_{m}^{2}$ with
\begin{align}
&  \mathcal{A}^{((d_{1},d_{2})\in D\oplus D\mapsto(d_{1},d_{2})\in D^{2}%
)}(x)\nonumber\\
&  =\mathcal{A}^{((d_{1},d_{2})\in D\oplus D\mapsto(d_{1},d_{2})\in D^{2})}(y)
\label{n2.3}%
\end{align}
By using the quasi-colimit diagram of small objects at page 92 of Lavendhomme
\cite{l1}, we are sure that there exists a unique $z\in\mathcal{A}_{m}%
^{D^{2}\oplus D}$ with
\begin{equation}
\mathcal{A}^{((d_{1},d_{2})\in D^{2}\mapsto(d_{1},d_{2},0)\in D^{2}\oplus
D)}(z)=x \label{n2.1}%
\end{equation}
and
\begin{equation}
\mathcal{A}^{((d_{1},d_{2})\in D^{2}\mapsto(d_{1},d_{2},d_{1}d_{2})\in
D^{2}\oplus D)}(z)=y \label{n2.2}%
\end{equation}
We define $y\overset{\cdot}{-}x\in\mathcal{A}_{m}^{1}$ to be $\mathcal{A}%
^{(d\in D\mapsto(0,0,d)\in D^{2}\oplus D)}(z)$.

\begin{proposition}
\label{tn2.2}Let $x,y\in\mathcal{A}^{2}$ abide by (\ref{n2.3}). Then we have
\begin{align*}
&  \mathcal{A}^{((d_{1},d_{2})\in D^{2}\mapsto(d_{2},d_{1})\in D^{2}%
)}(y)\overset{\cdot}{-}\mathcal{A}^{((d_{1},d_{2})\in D^{2}\mapsto(d_{2}%
,d_{1})\in D^{2})}(x)\\
&  =y\overset{\cdot}{-}x
\end{align*}

\end{proposition}

\begin{proof}
Let $z\in\mathcal{A}_{m}^{D^{2}\oplus D}$ obedient to (\ref{n2.1}) and
(\ref{n2.2}). Then we have
\begin{align*}
&  \mathcal{A}^{((d_{1},d_{2})\in D^{2}\mapsto(d_{1},d_{2},0)\in D^{2}\oplus
D)}\circ\mathcal{A}^{((d_{1},d_{2},d_{3})\in D^{2}\oplus D\mapsto(d_{2}%
,d_{1},d_{3})\in D^{2}\oplus D)}(z)\\
&  =\mathcal{A}^{((d_{1},d_{2})\in D^{2}\mapsto(d_{2},d_{1},0)\in D^{2}\oplus
D)}(z)\\
&  =\mathcal{A}^{((d_{1},d_{2})\in D^{2}\mapsto(d_{2},d_{1})\in D^{2})}%
\circ\mathcal{A}^{((d_{1},d_{2})\in D^{2}\mapsto(d_{1},d_{2},0)\in D^{2}\oplus
D)}(z)\\
&  =\mathcal{A}^{((d_{1},d_{2})\in D^{2}\mapsto(d_{2},d_{1})\in D^{2})}(x)
\end{align*}
while we have
\begin{align*}
&  \mathcal{A}^{((d_{1},d_{2})\in D^{2}\mapsto(d_{1},d_{2},d_{1}d_{2})\in
D^{2}\oplus D)}\circ\mathcal{A}^{((d_{1},d_{2},d_{3})\in D^{2}\oplus
D\mapsto(d_{2},d_{1},d_{3})\in D^{2}\oplus D)}(z)\\
&  =\mathcal{A}^{((d_{1},d_{2})\in D^{2}\mapsto(d_{2},d_{1},d_{1}d_{2})\in
D^{2}\oplus D)}(z)\\
&  =\mathcal{A}^{((d_{1},d_{2})\in D^{2}\mapsto(d_{2},d_{1})\in D^{2})}%
\circ\mathcal{A}^{((d_{1},d_{2})\in D^{2}\mapsto(d_{1},d_{2},d_{1}d_{2})\in
D^{2}\oplus D)}(z)\\
&  =\mathcal{A}^{((d_{1},d_{2})\in D^{2}\mapsto(d_{2},d_{1})\in D^{2})}(y)
\end{align*}
Therefore we have
\begin{align*}
&  \mathcal{A}^{((d_{1},d_{2})\in D^{2}\mapsto(d_{2},d_{1})\in D^{2}%
)}(y)\overset{\cdot}{-}\mathcal{A}^{((d_{1},d_{2})\in D^{2}\mapsto(d_{2}%
,d_{1})\in D^{2})}(x)\\
&  =\mathcal{A}^{(d\in D\mapsto(0,0,d)\in D^{2}\oplus D)}\circ\mathcal{A}%
^{((d_{1},d_{2},d_{3})\in D^{2}\oplus D\mapsto(d_{2},d_{1},d_{3})\in
D^{2}\oplus D)}(z)\\
&  =\mathcal{A}^{(d\in D\mapsto(0,0,d)\in D^{2}\oplus D)}(z)\\
&  =y\overset{\cdot}{-}x
\end{align*}
This completes the proof.
\end{proof}

\begin{proposition}
Let $x,y\in\mathcal{A}^{2}$ abide by (\ref{n2.3}). Then we have
\[
x\overset{\cdot}{-}y=-(y\overset{\cdot}{-}x)
\]

\end{proposition}

\begin{proof}
Let $z\in\mathcal{A}^{D^{2}\oplus D}$ abide by the conditions (\ref{n2.1}) and
(\ref{n2.2}). Let $u\in\mathcal{A}^{D^{2}\oplus D}$ be
\[
u=\mathcal{A}^{((d_{1},d_{2},d_{3})\in D^{2}\oplus D\mapsto(d_{1},d_{2}%
,d_{1}d_{2}-d_{3})\in D^{2}\oplus D)}(z)
\]
Then we have
\begin{align*}
&  \mathcal{A}^{((d_{1},d_{2})\in D^{2}\mapsto(d_{1},d_{2},0)\in D^{2}\oplus
D)}(u)\\
&  =\mathcal{A}^{((d_{1},d_{2})\in D^{2}\mapsto(d_{1},d_{2},0)\in D^{2}\oplus
D)}\circ\mathcal{A}^{((d_{1},d_{2},d_{3})\in D^{2}\oplus D\mapsto(d_{1}%
,d_{2},d_{1}d_{2}-d_{3})\in D^{2}\oplus D)}(z)\\
&  =\mathcal{A}^{((d_{1},d_{2})\in D^{2}\mapsto(d_{1},d_{2},d_{1}d_{2})\in
D^{2}\oplus D)}(z)\\
&  =y
\end{align*}
while we have
\begin{align*}
&  \mathcal{A}^{((d_{1},d_{2})\in D^{2}\mapsto(d_{1},d_{2},d_{1}d_{2})\in
D^{2}\oplus D)}(u)\\
&  =\mathcal{A}^{((d_{1},d_{2})\in D^{2}\mapsto(d_{1},d_{2},d_{1}d_{2})\in
D^{2}\oplus D)}\circ\mathcal{A}^{((d_{1},d_{2},d_{3})\in D^{2}\oplus
D\mapsto(d_{1},d_{2},d_{1}d_{2}-d_{3})\in D^{2}\oplus D)}(z)\\
&  =\mathcal{A}^{((d_{1},d_{2})\in D^{2}\mapsto(d_{1},d_{2},0)\in D^{2}\oplus
D)}(z)\\
&  =x
\end{align*}
Therefore we have
\begin{align*}
&  x\overset{\cdot}{-}y\\
&  =\mathcal{A}^{(d\in D\mapsto(0,0,d)\in D^{2}\oplus D)}(u)\\
&  =\mathcal{A}^{(d\in D\mapsto(0,0,d)\in D^{2}\oplus D)}\circ\mathcal{A}%
^{((d_{1},d_{2},d_{3})\in D^{2}\oplus D\mapsto(d_{1},d_{2},d_{1}d_{2}%
-d_{3})\in D^{2}\oplus D)}(z)\\
&  =\mathcal{A}^{(d\in D\mapsto(0,0,-d)\in D^{2}\oplus D)}(z)\\
&  =-(y\overset{\cdot}{-}x)
\end{align*}
This completes the proof.
\end{proof}

\begin{proposition}
Let $x,y\in\mathcal{A}^{2}$ abide by (\ref{n2.3}) with $a\in\mathbb{R}$. Then
we have
\[
a\underset{i}{\cdot}y\overset{\cdot}{-}a\underset{i}{\cdot}x=a(y\overset
{\cdot}{-}x)\text{ \ \ \ \ \ \ (}i=1,2\text{)}%
\]

\end{proposition}

\begin{proof}
Here we deal only with the case $i=1$, leaving the other case to the reader.
Let $z\in\mathcal{A}^{D^{2}\oplus D}$ abide by the conditions (\ref{n2.1}) and
(\ref{n2.2}). Let $u\in\mathcal{A}^{D^{2}\oplus D}$ be
\[
u=\mathcal{A}^{((d_{1},d_{2},d_{3})\in D^{2}\oplus D\mapsto(ad_{1}%
,d_{2},ad_{3})\in D^{2}\oplus D)}(z)
\]
Then we have
\begin{align*}
&  \mathcal{A}^{((d_{1},d_{2})\in D^{2}\mapsto(d_{1},d_{2},0)\in D^{2}\oplus
D)}(u)\\
&  =\mathcal{A}^{((d_{1},d_{2})\in D^{2}\mapsto(d_{1},d_{2},0)\in D^{2}\oplus
D)}\circ\mathcal{A}^{((d_{1},d_{2},d_{3})\in D^{2}\oplus D\mapsto(ad_{1}%
,d_{2},ad_{3})\in D^{2}\oplus D)}(z)\\
&  =\mathcal{A}^{((d_{1},d_{2})\in D^{2}\mapsto(ad_{1},d_{2},0)\in D^{2}\oplus
D)}(z)\\
&  =\mathcal{A}^{((d_{1},d_{2})\in D^{2}\mapsto(ad_{1},d_{2})\in D^{2})}%
\circ\mathcal{A}^{((d_{1},d_{2})\in D^{2}\mapsto(d_{1},d_{2},0)\in D^{2}\oplus
D)}(z)\\
&  =a\underset{1}{\cdot}x
\end{align*}
while we have
\begin{align*}
&  \mathcal{A}^{((d_{1},d_{2})\in D^{2}\mapsto(d_{1},d_{2},d_{1}d_{2})\in
D^{2}\oplus D)}(u)\\
&  =\mathcal{A}^{((d_{1},d_{2})\in D^{2}\mapsto(d_{1},d_{2},d_{1}d_{2})\in
D^{2}\oplus D)}\circ\mathcal{A}^{((d_{1},d_{2},d_{3})\in D^{2}\oplus
D\mapsto(ad_{1},d_{2},ad_{3})\in D^{2}\oplus D)}(z)\\
&  =\mathcal{A}^{((d_{1},d_{2})\in D^{2}\mapsto(ad_{1},d_{2},ad_{1}d_{2})\in
D^{2}\oplus D)}(z)\\
&  =\mathcal{A}^{((d_{1},d_{2})\in D^{2}\mapsto(ad_{1},d_{2})\in D^{2})}%
\circ\mathcal{A}^{((d_{1},d_{2})\in D^{2}\mapsto(d_{1},d_{2},d_{1}d_{2})\in
D^{2}\oplus D)}(z)\\
&  =a\underset{1}{\cdot}y
\end{align*}
Therefore we have
\begin{align*}
&  a\underset{i}{\cdot}y\overset{\cdot}{-}a\underset{i}{\cdot}x\\
&  =\mathcal{A}^{(d\in D\mapsto(0,0,d)\in D^{2}\oplus D)}(u)\\
&  =\mathcal{A}^{(d\in D\mapsto(0,0,d)\in D^{2}\oplus D)}\circ\mathcal{A}%
^{((d_{1},d_{2},d_{3})\in D^{2}\oplus D\mapsto(ad_{1},d_{2},ad_{3})\in
D^{2}\oplus D)}(z)\\
&  =\mathcal{A}^{(d\in D\mapsto(0,0,ad)\in D^{2}\oplus D)}(z)\\
&  =\mathcal{A}^{(d\in D\mapsto ad\in D)}\circ\mathcal{A}^{(d\in
D\mapsto(0,0,d)\in D^{2}\oplus D)}(z)\\
&  =a(y\overset{\cdot}{-}x)
\end{align*}
This completes the proof.
\end{proof}

\begin{lemma}
The following diagram is a quasi-colimit diagram:
\[%
\begin{array}
[c]{cccccccc}
&  & D^{2} &  & \overset{\mathbf{i}}{\leftarrow} & D\oplus D &  & \\
& \overset{\mathbf{i}}{\nearrow} &  & \overset{\mathbf{\varphi}_{1}}{\searrow}
&  &  & \overset{\mathbf{i}}{\searrow} & \\
D\oplus D &  &  &  & D^{2}\oplus D\oplus D & \overset{\mathbf{\varphi}_{3}%
}{\longleftarrow} &  & D^{2}\\
& \underset{\mathbf{i}}{\searrow} &  & \underset{\mathbf{\varphi}_{2}%
}{\nearrow} &  &  & \underset{\mathbf{i}}{\nearrow} & \\
&  & D^{2} &  & \underset{\mathbf{i}}{\leftarrow} & D\oplus D &  &
\end{array}
\]
where $\mathbf{i}:D\oplus D\rightarrow D^{2}$ is the canonical injection, and
$\mathbf{\varphi}_{1},\mathbf{\varphi}_{2},\mathbf{\varphi}_{3}:D^{2}%
\rightarrow D^{2}\oplus D\oplus D$ are defined to be
\begin{align*}
\mathbf{\varphi}_{1}(d_{1},d_{2})  &  =(d_{1},d_{2},0,0)\\
\mathbf{\varphi}_{2}(d_{1},d_{2})  &  =(d_{1},d_{2},d_{1}d_{2},0)\\
\mathbf{\varphi}_{3}(d_{1},d_{2})  &  =(d_{1},d_{2},0,d_{1}d_{2})
\end{align*}

\end{lemma}

\begin{proposition}
Let $x,y,z\in\mathcal{A}^{2}$ with
\begin{align*}
&  \mathcal{A}^{((d_{1},d_{2})\in D\oplus D\mapsto(d_{1},d_{2})\in D^{2}%
)}(x)\\
&  =\mathcal{A}^{((d_{1},d_{2})\in D\oplus D\mapsto(d_{1},d_{2})\in D^{2}%
)}(y)\\
&  =\mathcal{A}^{((d_{1},d_{2})\in D\oplus D\mapsto(d_{1},d_{2})\in D^{2})}(z)
\end{align*}
Then we have
\[
(y\overset{\cdot}{-}x)+(z\overset{\cdot}{-}y)+(x\overset{\cdot}{-}%
z)=\mathbf{0}%
\]

\end{proposition}

\begin{proof}
Let $u\in\mathcal{A}^{D^{2}\oplus D\oplus D}$ be the unique one such that
\begin{align*}
x  &  =\mathcal{A}^{((d_{1},d_{2})\in D^{2}\mapsto(d_{1},d_{2},0,0)\in
D^{2}\oplus D\oplus D)}(u)\\
y  &  =\mathcal{A}^{((d_{1},d_{2})\in D^{2}\mapsto(d_{1},d_{2},d_{1}%
d_{2},0)\in D^{2}\oplus D\oplus D)}(u)\\
z  &  =\mathcal{A}^{((d_{1},d_{2})\in D^{2}\mapsto(d_{1},d_{2},0,d_{1}%
d_{2})\in D^{2}\oplus D\oplus D)}(u)
\end{align*}
The unique existence of such $u\in\mathcal{A}^{D^{2}\oplus D\oplus D}$ is
guaranteed by the above lemma. Since we have
\begin{align*}
&  x\\
&  =\mathcal{A}^{((d_{1},d_{2})\in D^{2}\mapsto(d_{1},d_{2},0,0)\in
D^{2}\oplus D\oplus D)}(u)\\
&  =\mathcal{A}^{((d_{1},d_{2})\in D^{2}\mapsto(d_{1},d_{2},0)\in D^{2}\oplus
D)}\circ\mathcal{A}^{((d_{1},d_{2},d_{3})\in D^{2}\oplus D\mapsto(d_{1}%
,d_{2},d_{3},0)\in D^{2}\oplus D\oplus D)}(u)
\end{align*}
and
\begin{align*}
&  y\\
&  =\mathcal{A}^{((d_{1},d_{2})\in D^{2}\mapsto(d_{1},d_{2},d_{1}d_{2},0)\in
D^{2}\oplus D\oplus D)}(u)\\
&  =\mathcal{A}^{((d_{1},d_{2})\in D^{2}\mapsto(d_{1},d_{2},d_{1}d_{2})\in
D^{2}\oplus D)}\circ\mathcal{A}^{((d_{1},d_{2},d_{3})\in D^{2}\oplus
D\mapsto(d_{1},d_{2},d_{3},0)\in D^{2}\oplus D\oplus D)}(u)
\end{align*}
we have
\begin{align}
&  y\overset{\cdot}{-}x\nonumber\\
&  =\mathcal{A}^{(d\in D\mapsto(0,0,d)\in D^{2}\oplus D)}\circ\mathcal{A}%
^{((d_{1},d_{2},d_{3})\in D^{2}\oplus D\mapsto(d_{1},d_{2},d_{3},0)\in
D^{2}\oplus D\oplus D)}(u)\nonumber\\
&  =\mathcal{A}^{(d\in D\mapsto(0,0,d,0)\in D^{2}\oplus D\oplus D)}(u)
\label{n2.10}%
\end{align}
Since we have
\begin{align*}
&  y\\
&  =\mathcal{A}^{((d_{1},d_{2})\in D^{2}\mapsto(d_{1},d_{2},d_{1}d_{2},0)\in
D^{2}\oplus D\oplus D)}(u)\\
&  =\mathcal{A}^{((d_{1},d_{2})\in D^{2}\mapsto(d_{1},d_{2},0)\in D^{2}\oplus
D)}\circ\mathcal{A}^{((d_{1},d_{2},d_{3})\in D^{2}\oplus D\mapsto(d_{1}%
,d_{2},d_{1}d_{2}-d_{3},d_{3})\in D^{2}\oplus D\oplus D)}(u)
\end{align*}
and
\begin{align*}
&  z\\
&  =\mathcal{A}^{((d_{1},d_{2})\in D^{2}\mapsto(d_{1},d_{2},0,d_{1}d_{2})\in
D^{2}\oplus D\oplus D)}(u)\\
&  =\mathcal{A}^{((d_{1},d_{2})\in D^{2}\mapsto(d_{1},d_{2},d_{1}d_{2})\in
D^{2}\oplus D)}\circ\mathcal{A}^{((d_{1},d_{2},d_{3})\in D^{2}\oplus
D\mapsto(d_{1},d_{2},d_{1}d_{2}-d_{3},d_{3})\in D^{2}\oplus D\oplus D)}(u)
\end{align*}
we have
\begin{align}
&  z\overset{\cdot}{-}y\nonumber\\
&  =\mathcal{A}^{(d\in D\mapsto(0,0,d)\in D^{2}\oplus D)}\circ\mathcal{A}%
^{((d_{1},d_{2},d_{3})\in D^{2}\oplus D\mapsto(d_{1},d_{2},d_{1}d_{2}%
-d_{3},d_{3})\in D^{2}\oplus D\oplus D)}(u)\nonumber\\
&  =\mathcal{A}^{(d\in D\mapsto(0,0,-d,d)\in D^{2}\oplus D\oplus D)}(u)
\label{n2.11}%
\end{align}
Since we have
\begin{align*}
&  z\\
&  =\mathcal{A}^{((d_{1},d_{2})\in D^{2}\mapsto(d_{1},d_{2},0,d_{1}d_{2})\in
D^{2}\oplus D\oplus D)}(u)\\
&  =\mathcal{A}^{((d_{1},d_{2})\in D^{2}\mapsto(d_{1},d_{2},0)\in D^{2}\oplus
D)}\circ\mathcal{A}^{((d_{1},d_{2},d_{3})\in D^{2}\oplus D\mapsto(d_{1}%
,d_{2},0,d_{1}d_{2}-d_{3})\in D^{2}\oplus D\oplus D)}(u)
\end{align*}
and
\begin{align*}
&  x\\
&  =\mathcal{A}^{((d_{1},d_{2})\in D^{2}\mapsto(d_{1},d_{2},0,0)\in
D^{2}\oplus D\oplus D)}(u)\\
&  =\mathcal{A}^{((d_{1},d_{2})\in D^{2}\mapsto(d_{1},d_{2},d_{1}d_{2})\in
D^{2}\oplus D)}\circ\mathcal{A}^{((d_{1},d_{2},d_{3})\in D^{2}\oplus
D\mapsto(d_{1},d_{2},0,d_{1}d_{2}-d_{3})\in D^{2}\oplus D\oplus D)}(u)
\end{align*}
we have
\begin{align}
&  x\overset{\cdot}{-}z\nonumber\\
&  =\mathcal{A}^{(d\in D\mapsto(0,0,d)\in D^{2}\oplus D)}\circ\mathcal{A}%
^{((d_{1},d_{2},d_{3})\in D^{2}\oplus D\mapsto(d_{1},d_{2},0,d_{1}d_{2}%
-d_{3})\in D^{2}\oplus D\oplus D)}(u)\nonumber\\
&  =\mathcal{A}^{(d\in D\mapsto(0,0,0,-d)\in D^{2}\oplus D\oplus D)}(u)
\label{n2.12}%
\end{align}
Since we have
\begin{align*}
&  y\overset{\cdot}{-}x\\
&  =\mathcal{A}^{(d\in D\mapsto(0,0,d,0)\in D^{2}\oplus D\oplus D)}(u)\text{
\ \ [by (\ref{n2.10})]}\\
&  =\mathcal{A}^{(d\in D\mapsto(d,0)\in D\oplus D)}\circ\mathcal{A}%
^{((d_{1},d_{2})\in D\oplus D\mapsto(0,0,d_{1}-d_{2},d_{2})\in D^{2}\oplus
D\oplus D)}(u)
\end{align*}
and
\begin{align*}
&  z\overset{\cdot}{-}y\\
&  =\mathcal{A}^{(d\in D\mapsto(0,0,-d,d)\in D^{2}\oplus D\oplus D)}(u)\text{
\ \ \ [by (\ref{n2.11})]}\\
&  =\mathcal{A}^{(d\in D\mapsto(0,d)\in D\oplus D)}\circ\mathcal{A}%
^{((d_{1},d_{2})\in D\oplus D\mapsto(0,0,d_{1}-d_{2},d_{2})\in D^{2}\oplus
D\oplus D)}(u)
\end{align*}
we have
\begin{align}
&  (y\overset{\cdot}{-}x)+(z\overset{\cdot}{-}y)\nonumber\\
&  =\mathcal{A}^{(d\in D\mapsto(d,d)\in D\oplus D)}\circ\mathcal{A}%
^{((d_{1},d_{2})\in D\oplus D\mapsto(0,0,d_{1}-d_{2},d_{2})\in D^{2}\oplus
D\oplus D)}(u)\nonumber\\
&  =\mathcal{A}^{(d\in D\mapsto(0,0,0,d)\in D^{2}\oplus D\oplus D)}(u)
\label{n2.13}%
\end{align}
Since we have
\begin{align*}
&  (y\overset{\cdot}{-}x)+(z\overset{\cdot}{-}y)\\
&  =\mathcal{A}^{(d\in D\mapsto(0,0,0,d)\in D^{2}\oplus D\oplus D)}(u)\text{
\ \ \ [by (\ref{n2.13})]}\\
&  =\mathcal{A}^{(d\in D\mapsto(d,0)\in D\oplus D)}\circ\mathcal{A}%
^{((d_{1},d_{2})\in D\oplus D\mapsto(0,0,0,d_{1}-d_{2})\in D^{2}\oplus D\oplus
D)}(u)
\end{align*}
and
\begin{align*}
&  x\overset{\cdot}{-}z\\
&  =\mathcal{A}^{(d\in D\mapsto(0,0,0,-d)\in D^{2}\oplus D\oplus D)}(u)\text{
\ \ \ \ [by (\ref{n2.12})]}\\
&  =\mathcal{A}^{(d\in D\mapsto(0,d)\in D\oplus D)}\circ\mathcal{A}%
^{((d_{1},d_{2})\in D\oplus D\mapsto(0,0,0,d_{1}-d_{2})\in D^{2}\oplus D\oplus
D)}(u)
\end{align*}
we have
\begin{align*}
&  \{(y\overset{\cdot}{-}x)+(z\overset{\cdot}{-}y)\}+(x\overset{\cdot}{-}z)\\
&  =\mathcal{A}^{(d\in D\mapsto(d,d)\in D\oplus D)}\circ\mathcal{A}%
^{((d_{1},d_{2})\in D\oplus D\mapsto(0,0,0,d_{1}-d_{2})\in D^{2}\oplus D\oplus
D)}(u)\\
&  =\mathcal{A}^{(d\in D\mapsto(0,0,0,0)\in D^{2}\oplus D\oplus D)}(u)\\
&  =\mathbf{0}%
\end{align*}
This completes the proof.
\end{proof}

Let $x,y\in\mathcal{A}_{m}^{3}$ with
\begin{align}
&  \mathcal{A}^{((d_{1},d_{2},d_{3})\in D\times(D\oplus D)\mapsto(d_{1}%
,d_{2},d_{3})\in D^{3})}(x)\nonumber\\
&  =\mathcal{A}^{((d_{1},d_{2},d_{3})\in D\times(D\oplus D)\mapsto(d_{1}%
,d_{2},d_{3})\in D^{3})}(y) \label{n2.4}%
\end{align}
. By using the first quasi-colimit diagram of small objects in Lemma 2.1 of
Nishimura \cite{n11}, we are sure that there exists a unique $z\in
\mathcal{A}_{m}^{D^{4}\{(2,4),(3,4)\}}$ with
\[
\mathcal{A}^{((d_{1},d_{2},d_{3})\in D^{3}\mapsto(d_{1},d_{2},d_{3},0)\in
D^{4}\{(2,4),(3,4)\})}(z)=x
\]
and
\[
\mathcal{A}^{((d_{1},d_{2},d_{3})\in D^{3}\mapsto(d_{1},d_{2},d_{3},d_{2}%
d_{3})\in D^{4}\{(2,4),(3,4)\})}(z)=y
\]
We define $y\underset{1}{\overset{\cdot}{-}}x\in\mathcal{A}_{m}^{2}$ to be
$\mathcal{A}^{((d_{1},d_{2})\in D^{2}\mapsto(d_{1},0,0,d_{2})\in
D^{4}\{(2,4),(3,4)\})}(z)$.

Let $x,y\in\mathcal{A}_{m}^{3}$ with
\begin{align}
&  \mathcal{A}^{((d_{1},d_{2},d_{3})\in D^{3}\{(1,3)\}\mapsto(d_{1}%
,d_{2},d_{3})\in D^{3})}(x)\nonumber\\
&  =\mathcal{A}^{((d_{1},d_{2},d_{3})\in D^{3}\{(1,3)\}\mapsto(d_{1}%
,d_{2},d_{3})\in D^{3})}(y) \label{n2.5}%
\end{align}
By using the second quasi-colimit diagram of small objects in Lemma 2.1 of
Nishimura \cite{n11}, we are sure that there exists a unique $z\in
\mathcal{A}_{m}^{D^{4}\{(1,4),(3,4)\}}$ with
\[
\mathcal{A}^{((d_{1},d_{2},d_{3})\in D^{3}\mapsto(d_{1},d_{2},d_{3},0)\in
D^{4}\{(1,4),(3,4)\})}(z)=x
\]
and
\[
\mathcal{A}^{((d_{1},d_{2},d_{3})\in D^{3}\mapsto(d_{1},d_{2},d_{3},d_{1}%
d_{3})\in D^{4}\{(1,4),(3,4)\})}(z)=y
\]
We define $y\underset{2}{\overset{\cdot}{-}}x\in\mathcal{A}_{m}^{2}$ to be
$\mathcal{A}^{((d_{1},d_{2})\in D^{2}\mapsto(0,d_{1},0,d_{2})\in
D^{4}\{(1,4),(3,4)\})}(z)$.

Let $x,y\in\mathcal{A}_{m}^{3}$ with
\begin{align}
&  \mathcal{A}^{((d_{1},d_{2},d_{3})\in(D\oplus D)\times D\mapsto(d_{1}%
,d_{2},d_{3})\in D^{3})}(x)\label{n2.6}\\
&  =\mathcal{A}^{((d_{1},d_{2},d_{3})\in(D\oplus D)\times D\mapsto(d_{1}%
,d_{2},d_{3})\in D^{3})}(y)\nonumber
\end{align}
By using the third quasi-colimit diagram of small objects in Lemma 2.1 of
Nishimura \cite{n11}, we are sure that there exists a unique $z\in
\mathcal{A}_{m}^{D^{4}\{(1,4),(2,4)\}}$ with
\[
\mathcal{A}^{((d_{1},d_{2},d_{3})\in D^{3}\mapsto(d_{1},d_{2},d_{3},0)\in
D^{4}\{(1,4),(2,4)\})}(z)=x
\]
and
\[
\mathcal{A}^{((d_{1},d_{2},d_{3})\in D^{3}\mapsto(d_{1},d_{2},d_{3},d_{1}%
d_{2})\in D^{4}\{(1,4),(2,4)\})}(z)=y
\]
We define $y\underset{3}{\overset{\cdot}{-}}x\in\mathcal{A}_{m}^{2}$ to be
$\mathcal{A}^{((d_{1},d_{2})\in D^{2}\mapsto(0,0,d_{1},d_{2})\in
D^{4}\{(1,4),(2,4)\})}(z)$.

\begin{proposition}
\label{tn2.3}Let $x,y\in\mathcal{A}_{m}^{3}$.

\begin{enumerate}
\item If they satisfy (\ref{n2.4}), then we have
\begin{align*}
&  y\underset{1}{\overset{\cdot}{-}}x\\
&  =\mathcal{A}^{((d_{1},d_{2},d_{3})\in D^{3}\mapsto(d_{2},d_{1},d_{3})\in
D^{3})}(y)\underset{2}{\overset{\cdot}{-}}\mathcal{A}^{((d_{1},d_{2},d_{3})\in
D^{3}\mapsto(d_{2},d_{1},d_{3})\in D^{3})}(x)\\
&  =\mathcal{A}^{((d_{1},d_{2},d_{3})\in D^{3}\mapsto(d_{3},d_{2},d_{1})\in
D^{3})}(y)\underset{3}{\overset{\cdot}{-}}\mathcal{A}^{((d_{1},d_{2},d_{3})\in
D^{3}\mapsto(d_{3},d_{2},d_{1})\in D^{3})}(x)\\
&  =\mathcal{A}^{((d_{1},d_{2},d_{3})\in D^{3}\mapsto(d_{1},d_{3},d_{2})\in
D^{3})}(y)\underset{1}{\overset{\cdot}{-}}\mathcal{A}^{((d_{1},d_{2},d_{3})\in
D^{3}\mapsto(d_{1},d_{3},d_{2})\in D^{3})}(x)
\end{align*}

\item If they satisfy (\ref{n2.5}), then we have
\begin{align*}
&  y\underset{2}{\overset{\cdot}{-}}x\\
&  =\mathcal{A}^{((d_{1},d_{2},d_{3})\in D^{3}\mapsto(d_{2},d_{1},d_{3})\in
D^{3})}(y)\underset{1}{\overset{\cdot}{-}}\mathcal{A}^{((d_{1},d_{2},d_{3})\in
D^{3}\mapsto(d_{2},d_{1},d_{3})\in D^{3})}(x)\\
&  =\mathcal{A}^{((d_{1},d_{2},d_{3})\in D^{3}\mapsto(d_{3},d_{2},d_{1})\in
D^{3})}(y)\underset{2}{\overset{\cdot}{-}}\mathcal{A}^{((d_{1},d_{2},d_{3})\in
D^{3}\mapsto(d_{3},d_{2},d_{1})\in D^{3})}(x)\\
&  =\mathcal{A}^{((d_{1},d_{2},d_{3})\in D^{3}\mapsto(d_{1},d_{3},d_{2})\in
D^{3})}(y)\underset{3}{\overset{\cdot}{-}}\mathcal{A}^{((d_{1},d_{2},d_{3})\in
D^{3}\mapsto(d_{1},d_{3},d_{2})\in D^{3})}(x)
\end{align*}

\item If they satisfy (\ref{n2.6}), then we have
\begin{align*}
&  y\underset{3}{\overset{\cdot}{-}}x\\
&  =\mathcal{A}^{((d_{1},d_{2},d_{3})\in D^{3}\mapsto(d_{2},d_{1},d_{3})\in
D^{3})}(y)\underset{3}{\overset{\cdot}{-}}\mathcal{A}^{((d_{1},d_{2},d_{3})\in
D^{3}\mapsto(d_{2},d_{1},d_{3})\in D^{3})}(x)\\
&  =\mathcal{A}^{((d_{1},d_{2},d_{3})\in D^{3}\mapsto(d_{3},d_{2},d_{1})\in
D^{3})}(y)\underset{1}{\overset{\cdot}{-}}\mathcal{A}^{((d_{1},d_{2},d_{3})\in
D^{3}\mapsto(d_{3},d_{2},d_{1})\in D^{3})}(x)\\
&  =\mathcal{A}^{((d_{1},d_{2},d_{3})\in D^{3}\mapsto(d_{1},d_{3},d_{2})\in
D^{3})}(y)\underset{2}{\overset{\cdot}{-}}\mathcal{A}^{((d_{1},d_{2},d_{3})\in
D^{3}\mapsto(d_{1},d_{3},d_{2})\in D^{3})}(x)
\end{align*}

\end{enumerate}
\end{proposition}

\begin{proof}
The proof is similar to that in Proposition \ref{tn2.2}. The details can
safely be left to the reader.
\end{proof}

Now we have

\begin{theorem}
\label{tn2.4}The four strong differences $\overset{\cdot}{-}$, $\underset
{1}{\overset{\cdot}{-}}$, $\underset{2}{\overset{\cdot}{-}}$ and $\underset
{3}{\overset{\cdot}{-}}$ satisfy the general Jacobi identity. I.e., given
$x_{123},x_{132},x_{213},x_{231},x_{312},x_{321}\in\mathcal{A}^{3}$, as long
as the following three expressions are well defined, they sum up only to
vanish:
\begin{align*}
&  (x_{123}\overset{\cdot}{\underset{1}{-}}x_{132})\overset{\cdot}{-}%
(x_{231}\overset{\cdot}{\underset{1}{-}}x_{321})\\
&  (x_{231}\overset{\cdot}{\underset{2}{-}}x_{213})\overset{\cdot}{-}%
(x_{312}\overset{\cdot}{\underset{2}{-}}x_{132})\\
&  (x_{312}\overset{\cdot}{\underset{3}{-}}x_{321})\overset{\cdot}{-}%
(x_{123}\overset{\cdot}{\underset{3}{-}}x_{213})
\end{align*}

\end{theorem}

\begin{proof}
The theorem was already proved in case of the standard Nishimura algebroid
$\mathcal{S}_{M}$ in Nishimura's \cite{n12}, \S 3. What we should do is only
to reformulate the above proof genuinely in terms of diagrams. The details can
safely be left to the reader.
\end{proof}

\subsection{\textit{Nishimura algebroids}$_{3}$}

\begin{definition}
A Nishimura algebroid$_{2}$\ $\mathcal{A}$\ over $M$ is called a
\textit{Nishimura algebroid}$_{3}$ \textit{over} $M$ providing that it is
endowed with a natural transformation $\mathbf{a}$ from $\mathcal{A}$ to the
standard Nishimura algebroid$_{2}$ $\mathcal{S}_{M}$ to be called the anchor
natural transformation.
\end{definition}

\begin{example}
The standard Nishimura algebroid$_{2}$ $\mathcal{S}_{M}$ over $M$ is
canonically a Nishimura algebroid$_{3}$ over $M$ endowed with the identity
natural transformation of $\mathcal{S}_{M}$.
\end{example}

\begin{example}
Let $G$ be a groupoid over $M$. Then the Nishimura algebroid$_{2}$
$\mathcal{A}G$ over $M$ is a Nishimura algebroid$_{3}$ over $M$ endowed with
the anchor natural transformation assigning $\mathbf{a}_{G}^{\mathcal{D}%
}:\mathcal{A}^{\mathcal{D}}\rightarrow M^{\mathcal{D}}$ to each object
$\mathcal{D}$ in $\mathbf{Simp}$.
\end{example}

\subsection{\textit{Nishimura Algebroids}$_{4}$}

We denote by $\otimes_{\mathcal{A}}$, or more simply by $\otimes$, the
contravariant functor which assigns $\mathcal{D}_{1}\otimes\mathcal{D}%
_{2}=\{(\zeta,x)\in(\mathcal{A}^{\mathcal{D}_{2}})^{\mathcal{D}_{1}}%
\times\mathcal{A}^{\mathcal{D}_{1}}\mid\mathbf{a}(x)=\pi^{\mathcal{D}_{1}%
}(\zeta)\}$ to each object $(\mathcal{D}_{1},\mathcal{D}_{2})$ in
$\mathbf{Simp}\times\mathbf{Simp}$ and which assigns $f\otimes g=(\zeta
\in(\mathcal{A}^{\mathcal{D}_{2}})^{\mathcal{D}_{1}}\mapsto\mathcal{A}%
^{f}\circ\zeta\circ\mathcal{A}^{g}\in(\mathcal{A}^{\mathcal{D}_{2}^{\prime}%
})^{\mathcal{D}_{1}^{\prime}},\mathcal{A}^{g}):\mathcal{D}_{1}\otimes
\mathcal{D}_{2}\rightarrow\mathcal{D}_{1}^{\prime}\otimes\mathcal{D}%
_{2}^{\prime}$ to each morphism $(f,g):(\mathcal{D}_{1}^{\prime}%
,\mathcal{D}_{2}^{\prime})\rightarrow(\mathcal{D}_{1},\mathcal{D}_{2})$ in
$\mathbf{Simp}\times\mathbf{Simp}$, where $(\mathcal{A}^{\mathcal{D}_{2}%
})^{\mathcal{D}_{1}}$ denotes the space of mappings from the infinitesimal
space $\mathcal{D}_{1}$ to $\mathcal{A}^{\mathcal{D}_{2}}$, and $\pi
^{\mathcal{D}_{1}}(\zeta)$ assigns $\pi(\zeta(d))$ to each $d\in
\mathcal{D}_{1}$. We denote by $\widetilde{\otimes}_{\mathcal{A}}$, or more
simply by $\widetilde{\otimes\text{,}}$ the contravariant functor which
assigns $\mathcal{D}_{1}\widetilde{\otimes}\mathcal{D}_{2}=\mathcal{A}%
^{\mathcal{D}_{1}\times\mathcal{D}_{2}}$ to each object $(\mathcal{D}%
_{1},\mathcal{D}_{2})$ in $\mathbf{Simp}\times\mathbf{Simp}$ and which assigns
$f\widetilde{\otimes}g=\mathcal{A}^{f\times g}:\mathcal{A}^{\mathcal{D}%
_{1}\times\mathcal{D}_{2}}\rightarrow\mathcal{A}^{\mathcal{D}_{1}^{\prime
}\times\mathcal{D}_{2}^{\prime}}$ to each morphism $(f,g):(\mathcal{D}%
_{1}^{\prime},\mathcal{D}_{2}^{\prime})\rightarrow(\mathcal{D}_{1}%
,\mathcal{D}_{2})$ in $\mathbf{Simp}\times\mathbf{Simp}$.

\begin{definition}
\label{d2.4}A Nishimura algebroid$_{3}$ $\mathcal{A}$ over $M$ is called a
\textit{Nishimura algebroid}$_{4}$ \textit{over} $M$ providing that it is
endowed with a natural isomorphism $\ast_{\mathcal{A}}$ (denoted more simply
$\ast$ unless there is possible confusion) from the contravariant functor
$\otimes$ to the contravariant functor $\widetilde{\otimes}$ abiding by the
following conditions:

\begin{enumerate}
\item For any $(\zeta,x)\in\mathcal{D}_{1}\otimes\mathcal{D}_{2}$ with
$(\mathcal{D}_{1},\mathcal{D}_{2})$ in $\mathbf{Simp}\times\mathbf{Simp}$, we
have
\[
\pi(\zeta\ast x)=\pi(x)
\]
and
\[
\mathbf{a}(\zeta\ast x)=\mathbf{a}^{\mathcal{D}_{1}}(\zeta)
\]
where $\mathbf{a}^{\mathcal{D}_{1}}(\zeta)$ assigns $\mathbf{a}(\zeta
(d_{1}))(d_{2})$ to each $(d_{1},d_{2})\in\mathcal{D}_{1}\times\mathcal{D}%
_{2}$.

\item Let $\mathbf{i}_{j}:\mathcal{D}_{j}\rightarrow\mathcal{D}_{1}%
\times\mathcal{D}_{2}$ be the canonical injection with $\mathbf{p}%
_{j}:\mathcal{D}_{1}\times\mathcal{D}_{2}\rightarrow\mathcal{D}_{j}$ the
canonical projection ($j=1,2$). Then we have
\[
\mathcal{A}^{\mathbf{i}_{1}}(\zeta\ast x)=x
\]
and
\[
\mathcal{A}^{\mathbf{i}_{2}}(\zeta\ast x)=\zeta(0_{\mathcal{D}_{2}})
\]
for any $(\zeta,x)\in\mathcal{D}_{1}\otimes\mathcal{D}_{2}$, while we have
\[
\mathcal{A}^{\mathbf{p}_{1}}(y)=(d\in\mathcal{D}_{1}\mapsto\mathbf{0}%
_{(\mathbf{a}y)(d)}^{\mathcal{D}_{2}})\ast y
\]
for any $y\in\mathcal{A}^{\mathcal{D}_{1}}$ and
\[
\mathcal{A}^{\mathbf{p}_{2}}(z)=(d\in\mathcal{D}_{1}\mapsto z)\ast
\mathbf{0}_{\pi(z)}^{\mathcal{D}_{1}}%
\]
for any $z\in\mathcal{A}^{\mathcal{D}_{2}}$.

\item Let $f\in\mathbb{R}^{\mathcal{D}}$. For any $(\zeta,x)\in(\mathcal{D}%
_{1}\times...\times\mathcal{D}_{n})\otimes\mathcal{D}$, we have
\begin{align*}
&  \mathcal{A}^{((d,d_{1},...,d_{n})\in\mathcal{D\times D}_{1}\times
...\times\mathcal{D}_{n}\mapsto(d,d_{1},...,d_{i-1},f(d)d_{i},d_{i+1}%
,...,d_{n})\in\mathcal{D\times D}_{1}\times...\times\mathcal{D}_{n})}%
(\zeta\ast x)\\
&  =\{d\in\mathcal{D}\mapsto\mathcal{A}^{((d_{1},...,d_{n})\in\mathcal{D}%
_{1}\times...\times\mathcal{D}_{n}\mapsto(d_{1},...,d_{i-1},f(d)d_{i}%
,d_{i+1},...,d_{n})\in\mathcal{D}_{1}\times...\times\mathcal{D}_{n})}%
\zeta(d)\in\mathcal{A}^{\mathcal{D}_{1}\times...\times\mathcal{D}_{n}}\}\ast
x\text{ \ \ \ }\\
\text{(}1  &  \leq i\leq n\text{)}%
\end{align*}

\item For any $x\in\mathcal{A}^{\mathcal{D}_{1}}$, any $\zeta_{1}%
\in(\mathcal{A}^{\mathcal{D}_{2}})^{\mathcal{D}_{1}}$ and any $\zeta_{2}%
\in(\mathcal{A}^{\mathcal{D}_{3}})^{\mathcal{D}_{1}\times\mathcal{D}_{2}}$
with $\mathbf{a}(x)=\pi^{\mathcal{D}_{1}}(\zeta_{1})$ and $\mathbf{a}%
^{\mathcal{D}_{1}}(\zeta_{1})=\pi^{\mathcal{D}_{1}\times\mathcal{D}_{2}}%
(\zeta_{2})$, we have
\[
\zeta_{2}\ast(\zeta_{1}\ast x)=(\zeta_{2}\ast^{D}\zeta_{1})\ast x
\]
where $\zeta_{2}\ast^{\mathcal{D}_{1}}\zeta_{1}\in(\mathcal{A}^{\mathcal{D}%
_{1}\times\mathcal{D}_{2}})^{\mathcal{D}_{1}}$ is defined to be
\[
(\zeta_{2}\ast^{\mathcal{D}_{1}}\zeta_{1})(d)=\zeta_{2}(d,\cdot)\ast\zeta
_{1}(d)
\]
for any $d\in\mathcal{D}_{1}$.
\end{enumerate}
\end{definition}

\begin{remark}
\label{t2.4.1}What we require in our definition of Nishimura algebroid$_{4}$
over $M$ is that while multiplication seen in groupoids is no longer in view
in Nishimura algebroids, the remnants of multiplication and its associativity
are to be still in view. Multiplication seems completely lost in the
traditional definition of Lie algebroid.
\end{remark}

\begin{example}
The standard Nishimura algebroid$_{3}$ $\mathcal{S}_{M}$ over $M$ is
canonically a Nishimura algebroid$_{4}$ over $M$ provided that $\zeta
\ast_{\mathcal{S}_{M}}x\in\mathcal{S}_{M}^{\mathcal{D}_{1}\times
\mathcal{D}_{2}}$ is defined to be
\[
(d_{1},d_{2})\in\mathcal{D}_{1}\times\mathcal{D}_{2}\mapsto\zeta(d_{1}%
)(d_{2})\in M
\]

\end{example}

\begin{example}
Let $G$ be a groupoid over $M$. The Nishimura algebroid$_{3}$ $\mathcal{A}G$
over $M$ is a Nishimura algebroid$_{4}$ over $M$ provided that $\zeta
\ast_{\mathcal{A}G}x\in(\mathcal{A}G)^{\mathcal{D}_{1}\times\mathcal{D}_{2}}$
is defined to be
\[
(d_{1},d_{2})\in\mathcal{D}_{1}\times\mathcal{D}_{2}\mapsto\zeta(d_{1}%
)(d_{2})x(d_{1})\in G
\]

\end{example}

Now we give some results holding for any Nishimura algebroid$_{4}$
$\mathcal{A}$ over $M$.

\begin{proposition}
\label{t2.4.2}There is a bijective correspondence between the mappings
$\Phi:D\rightarrow\mathcal{A}_{m}^{1}$ and the elements $x\in\mathcal{A}%
_{m}^{2} $ with $\mathcal{A}^{(d\in D\mapsto(d,0)\in D^{2})}(x)=\mathbf{0}%
_{m}^{D}$.
\end{proposition}

\begin{proof}
This follows simply from the first condition in the definition of Nishimura
algebroid$_{4}$ over $M$, which claims that the assignment of $\Phi
\ast\mathbf{0}_{m}^{D}\in\mathcal{A}_{m}^{2}$ to each mapping $\Phi
:D\rightarrow\mathcal{A}_{m}^{1}$ gives such a bijective correspondence.
\end{proof}

It is easy to see that

\begin{lemma}
\label{l2.1}Let $\mathbf{p}_{1}:\mathcal{D}_{1}\times\mathcal{D}%
_{2}\rightarrow\mathcal{D}_{1}$ be the canonical projection as in the second
condition of Definition \ref{d2.4}. Then we have
\[
\mathcal{A}^{\mathbf{p}_{1}}(\mathbf{0}_{m}^{\mathcal{D}_{1}})=\mathbf{0}%
_{m}^{\mathcal{D}_{1}\times\mathcal{D}_{2}}%
\]

\end{lemma}

As an easy consequence of the above proposition, we have

\begin{theorem}
\label{t2.4.3}Given a Nishimura algebroid$_{4}$ $\mathcal{A}$ over $M$ with
$m\in M$, the $\mathbb{R}$-module $\mathcal{A}_{m}^{1}$ is Euclidean.
\end{theorem}

\begin{proof}
We have already proved that $\mathcal{A}_{m}^{1}$ is naturally an $\mathbb{R}
$-module. Let $\varphi:D\rightarrow\mathcal{A}_{m}^{1}$ be a mapping. We will
consider another mapping $\Phi:D\rightarrow\mathcal{A}_{m}^{1}$ defined to be
\[
\Phi(d)=\varphi(d)-\varphi(0)
\]
for any $d\in D$. Let us consider $x=\Phi\ast\mathbf{0}_{m}^{D}\in
\mathcal{A}_{m}^{2}$. We have $\mathcal{A}^{(d\in D\mapsto(d,0))}%
(x)=\mathbf{0}_{m}^{D}$, while it is easy to see that $\mathcal{A}^{(d\in
D\mapsto(0,d))}(x)=\Phi(0)=\mathbf{0}_{m}^{D}$. Therefore there is a unique
$y\in\mathcal{A}_{m}^{1}$ with $\mathcal{A}^{((d_{1},d_{2})\in D^{2}\mapsto
d_{1}d_{2}\in D)}(y)=x$. Let us consider $\mathcal{A}^{((d_{1},d_{2})\in
D^{2}\mapsto d_{2}\in D)}(y)=y\ast\mathbf{0}_{m}^{D}\in\mathcal{A}_{m}^{2}$.
Then it is easy to see that
\begin{align*}
(d  &  \in D\longmapsto dy)\ast\mathbf{0}_{m}^{D}\\
&  =\mathcal{A}^{((d_{1},d_{2})\in D^{2}\longmapsto(d_{1},d_{1}d_{2})\in
D^{2})}(\mathcal{A}^{((d_{1},d_{2})\in D^{2}\mapsto d_{2}\in D)}(y))\\
&  =\mathcal{A}^{((d_{1},d_{2})\in D^{2}\mapsto d_{1}d_{2}\in D)}(y)\\
&  =x
\end{align*}
Therefore we have $\Phi\ast\mathbf{0}_{m}^{D}=(d\in D\longmapsto
dy)\ast\mathbf{0}_{m}^{D}$, which implies that
\[
\varphi(d)-\varphi(0)=dy
\]
for any $d\in D$. To see the uniqueness of such $y\in\mathcal{A}_{m}^{1}$, let
us suppose that some $z\in\mathcal{A}_{m}^{1}$ satisfies
\[
dz=\mathbf{0}_{m}^{D}%
\]
for any $d\in D$. Since $z\ast\mathbf{0}_{m}^{D}=\mathcal{A}^{((d_{1}%
,d_{2})\in D^{2}\mapsto d_{2}\in D)}(z)$, we have
\begin{align*}
&  (d\in D\rightarrow\mathbf{0}_{m}^{D})\ast\mathbf{0}_{m}^{D}\\
&  =(d\in D\longmapsto dz)\ast\mathbf{0}_{m}^{D}\\
&  =\mathcal{A}^{((d_{1},d_{2})\in D^{2}\longmapsto(d_{1},d_{1}d_{2})\in
D^{2})}(\mathcal{A}^{((d_{1},d_{2})\in D^{2}\mapsto d_{2}\in D)}(z))\\
&  =\mathcal{A}^{((d_{1},d_{2})\in D^{2}\mapsto d_{1}d_{2}\in D)}(z)
\end{align*}
Since $(d\in D\rightarrow\mathbf{0}_{m}^{D})\ast\mathbf{0}_{m}^{D}%
=\mathbf{0}_{m}^{D^{2}}$ by Lemma \ref{l2.1} and the second condition of
Definition \ref{d2.4}, the desired uniqueness follows from Proposition 1
(\S 2.2) of Lavendhomme \cite{l1}.
\end{proof}

Now we will discuss the relationship between $\ast$\ and strong differences.

\begin{proposition}
\label{t2.6}

\begin{enumerate}
\item For any $\zeta_{1},\zeta_{2}\in(\mathcal{A}^{2})^{D}$ and any
$x\in\mathcal{A}^{1}$ with
\[
\mathbf{a}(x)=\pi^{D}(\zeta_{1})=\pi^{D}(\zeta_{2})
\]
and
\begin{align*}
&  \mathcal{A}^{((d_{1},d_{2})\in D\oplus D\mapsto(d_{1},d_{2})\in D^{2}%
)}(\zeta_{1}(d))\\
&  =\mathcal{A}^{((d_{1},d_{2})\in D\oplus D\mapsto(d_{1},d_{2})\in D^{2}%
)}(\zeta_{2}(d))
\end{align*}
for any $d\in D$, we have
\[
(\zeta_{2}\overset{\cdot}{-}\zeta_{1})\ast x=\zeta_{2}\ast x\underset
{1}{\overset{\cdot}{-}}\zeta_{1}\ast x
\]
where $\zeta_{2}\overset{\cdot}{-}\zeta_{1}\in(\mathcal{A}^{1})^{D}$ is
defined to be
\[
(\zeta_{2}\overset{\cdot}{-}\zeta_{1})(d)=\zeta_{2}(d)\overset{\cdot}{-}%
\zeta_{1}(d)
\]
for any $d\in D$.

\item For any $x,y\in\mathcal{A}^{2}$ and any $\zeta\in(\mathcal{A}%
^{1})^{D^{2}\oplus D}$ with
\[
\mathbf{a}(x)=(d_{1},d_{2})\in D^{2}\mapsto\pi(\zeta(d_{1},d_{2},0))
\]
\[
\mathbf{a}(y)=(d_{1},d_{2})\in D^{2}\mapsto\pi(x(d_{1},d_{2},d_{1}d_{2}))
\]
and
\begin{align*}
&  \mathcal{A}^{((d_{1},d_{2})\in D\oplus D\mapsto(d_{1},d_{2})\in D^{2}%
)}(x)\\
&  =\mathcal{A}^{((d_{1},d_{2})\in D\oplus D\mapsto(d_{1},d_{2})\in D^{2})}(y)
\end{align*}
we have
\begin{align*}
&  \mathcal{A}^{((d_{1},d_{2})\in D\oplus D\mapsto(d_{1},d_{2})\in D^{2})}\\
(\{\zeta\circ(d  &  \in D\mapsto(0,0,d)\in D^{2}\oplus D)\}\ast(y\overset
{\cdot}{-}x))\\
&  =\{\zeta\circ((d_{1},d_{2})\in D^{2}\mapsto(d_{1},d_{2},d_{1}d_{2})\in
D^{2}\oplus D)\}\ast y\underset{3}{\overset{\cdot}{-}}\\
\{\zeta\circ((d_{1},d_{2})  &  \in D^{2}\mapsto(d_{1},d_{2},0)\in D^{2}\oplus
D)\}\ast x
\end{align*}

\end{enumerate}
\end{proposition}

\begin{proof}
It suffices to note that given an object $\mathcal{D}$ in $\mathbf{Simp}$, the
contravariant functor $\widetilde{\otimes}\mathcal{D}$ (resp. $\mathcal{D}%
\widetilde{\otimes}$) and therefore the functor $\otimes\mathcal{D}$ (resp.
$\mathcal{D}\otimes$) map every quasi-colimit diagram of small objects in
$\mathbf{Simp}$ to a limit diagram. Therefore the proof is merely a
reformulation of Proposition 2.6 of Nishimura \cite{n11}. The details can
safely be left to the reader.
\end{proof}

\subsection{\textit{Nishimura Algebroids}$_{5}$}

\begin{definition}
A Nishimura algebroid$_{4}$ $\mathcal{A}$ over $M$ is called a
\textit{Nishimura algebroid}$_{5}$ \textit{over} $M$ providing that the anchor
natural transformation $\mathbf{a}$ from $\mathcal{A}$ to the standard
Nishimura algebroid$_{4}$ $\mathcal{S}_{M}$ is a homomorphism of Nishimura
algebroids$_{4}$ over $M$. In other words, a Nishimura algebroid$_{4}$
$\mathcal{A}$ over $M$ is a \textit{Nishimura algebroid}$_{5}$ \textit{over}
$M$ providing that for any $(\zeta,x)\in\mathcal{D}_{1}\otimes_{\mathcal{A}%
}\mathcal{D}_{2}$ with $(\mathcal{D}_{1},\mathcal{D}_{2})$ in $\mathbf{Simp}%
\times\mathbf{Simp}$, we have
\[
\mathbf{a}(\zeta\ast_{\mathcal{A}}x)=\mathbf{a}^{\mathcal{D}_{1}}(\zeta
)\ast_{\mathcal{S}_{M}}\mathbf{a}(x)
\]

\end{definition}

\begin{example}
It is trivial to see that the standard Nishimura algebroid$_{4}$
$\mathcal{S}_{M}$ over $M$ is a \textit{Nishimura algebroid}$_{5}$
\textit{over} $M$, since $\mathbf{a}$ is the identity transformation.
\end{example}

\begin{example}
Let $G$ be a groupoid over $M$. It is easy to see that the Nishimura
algebroid$_{4}$ $\mathcal{A}G$ over $M$\ is a Nishimura algebroid$_{5}$ over
$M$. It is also easy to see that a homomorphism $\varphi:G\rightarrow
G^{\prime}$ of groupoids over $M$\ naturally gives rise to a homomorphism
$\mathcal{A}\varphi:\mathcal{A}G\rightarrow\mathcal{A}G^{\prime}$of Nishimura
algebroids$_{5}$ over $M$. Thus we obtain a functor $\mathcal{A}$\ from the
category of groupoids over $M$ to the category of Nishimura algebroids$_{5}$
over $M$.
\end{example}

The following proposition should be obvious.

\begin{proposition}
Let $\varphi:\mathcal{A}\rightarrow\mathcal{A}^{\prime}$ be a homomorphism of
Nishimura algebroids$_{5}$ over $M$. Then its kernel at each $m\in M$, denoted
by $\mathrm{\ker}_{m}\varphi$, assigning $(\mathrm{\ker}_{m}\varphi
)^{\mathcal{D}}=\{x\in\mathcal{A}^{\mathcal{D}}\mid\varphi(x)=0_{m}%
^{\mathcal{D}}\}$ to each object $\mathcal{D}$ in $\mathbf{Simp}$ and
assigning the restriction $(\mathrm{\ker}_{m}\varphi)^{f}:(\mathrm{\ker}%
_{m}\varphi)^{\mathcal{D}^{\prime}}\rightarrow(\mathrm{\ker}_{m}%
\varphi)^{\mathcal{D}}$ of $\mathcal{A}^{f}:\mathcal{A}^{\mathcal{D}^{\prime}%
}\rightarrow\mathcal{A}^{\mathcal{D}}$ to each morphism
$f:\mathcal{D\rightarrow D}^{\prime}$ in $\mathbf{Simp}$ is naturally a
Nishimura algebroid$_{5}$ over a single point.
\end{proposition}

\subsection{\textit{Nishimura Algebroids}$_{6}$}

Let $\mathcal{A}$ be a Nishimura algebroid$_{5}$ over $M$. Since the anchor
natural transformation $\mathbf{a}_{\mathcal{A}}:\mathcal{A}\rightarrow
\mathcal{S}_{M}$ is really a homomorphism of Nishimura algebroids$_{5}$ over
$M$, its kernel $\mathrm{\ker}_{m}\mathbf{a}_{\mathcal{A}}$ at each $m\in M$
is a Nishimura algebroid$_{5}$ over a single point by dint of the last
proposition of the previous subsection. By collecting $\mathrm{\ker}%
_{m}\mathbf{a}_{\mathcal{A}}$ over all $m\in M$, we obtain a bundle of
Nishimura algebroids$_{5}$ over a single point, which is called the
\textit{inner subalgebroid} of $\mathcal{A}$ and which is denoted by
$\mathbf{I}\mathcal{A} $. The reader should note that the inner subalgebroid
$\mathbf{I}\mathcal{A}$ of $\mathcal{A}$ can naturally be reckoned as a
Nishimura algebroid$_{5}$ over $M$ (as a subalgebroid of $\mathcal{A}$ in a
natural sense). In the next definition we will consider the frame groupoid of
Nishimura algebroids$_{5}$ over a single point for $\mathbf{I}\mathcal{A}$,
which is denoted by $\Phi_{Nishi_{5}}(\mathbf{I}\mathcal{A})$.

\begin{definition}
A Nishimura algebroid$_{5}$ $\mathcal{A}$ over $M$ is called a
\textit{Nishimura algebroid}$_{6}$ \textit{over} $M$ providing that it is
endowed with a homomorphism $\mathrm{ad}_{\mathcal{A}}$ (usually written
simply $\mathrm{ad}$) of Nishimura algebroids$_{5}$ over $M$ from
$\mathcal{A}$ to $\mathcal{A}(\Phi_{Nishi_{5}}(\mathbf{I}\mathcal{A}))$
abiding by the following condition:

\begin{enumerate}
\item We have
\[
\mathrm{ad}(x)(d_{1})\circ\mathrm{ad}(y)(d_{2})=(\mathrm{ad}((\mathrm{ad}%
(x)(d_{1}))(y)))(d_{2})\circ\mathrm{ad}(x)(d_{1})
\]
for any objects $\mathcal{D}_{1},\mathcal{D}_{2}$ in $\mathbf{Simp}$, any
$d_{1}\in\mathcal{D}_{1}$, any $d_{2}\in\mathcal{D}_{2}$, any $x\in
\mathcal{A}^{\mathcal{D}_{1}}$ and any $y\in(\mathbf{I}\mathcal{A}%
)^{\mathcal{D}_{2}}$ with $\pi(x)=\pi(y)$.

\item Given $x,y\in(\mathbf{I}\mathcal{A)}^{1}$ with $\pi(x)=\pi(y)$, we have
\[
(\mathrm{ad}(x))(d)(y)-y=d[x,y]
\]
for any $d\in D$.
\end{enumerate}
\end{definition}

\begin{example}
Since the inner subalgebroid $\mathbf{I}\mathcal{S}_{M}$ of the standard
Nishimura algebroid$_{5}$ $\mathcal{S}_{M}$ is trivial, $\mathcal{S}_{M}$ is
trivially a \textit{Nishimura algebroid}$_{6}$ \textit{over} $M$.
\end{example}

\begin{example}
Let $G$ be a groupoid over $M$. By assigning a mapping
\[
y\in(\mathbf{I}G)_{\alpha x}\mapsto xyx^{-1}\in(\mathbf{I}G)_{\beta x}%
\]
to each $x\in G$, we get a homomorphism of groupoids over $M$ from $G$ to
$\Phi_{grp}(\mathbf{I}G)$, which naturally gives rise to a homomorphism of
groupoids over $M$ from $G$ to $\Phi_{Nishi_{5}}(\mathcal{A}(\mathbf{I}G))$.
Since $\mathcal{A}(\mathbf{I}G)$ and $\mathbf{I}(\mathcal{A}G)$ can naturally
be identified, we have a homomorphism of groupoids over $M$ from $G $ to
$\Phi_{Nishi_{5}}(\mathbf{I}(\mathcal{A}G))$, to which we apply the functor
$\mathcal{A}$ so as to get the desired $\mathrm{ad}_{\mathcal{A}G}$ as a
homomorphism of Nishimura algebroids$_{5}$ over $M$ from $\mathcal{A}G$ to
$\mathcal{A}(\Phi_{Nishi_{5}}(\mathbf{I(}\mathcal{A}G\mathcal{)}))$.
\end{example}

\section{Totally Intransitive Nishimura Algebroids}

\begin{definition}
A Nishimura algebroid $\mathcal{A}$ over $M$ is said to be totally
intransitive providing that its anchor natural transformation $\mathbf{a}%
_{\mathcal{A}}$ is trivial, i.e.,
\[
\mathbf{a}_{\mathcal{A}}(x)=\mathbf{0}_{m}^{\mathcal{D}}%
\]
for any $m\in M$, any object $\mathcal{D}$ in $\mathbf{Simp}$ and any
$x\in\mathcal{A}_{m}^{\mathcal{D}}$.
\end{definition}

\begin{remark}
A totally intransitive Nishimura algebroid $\mathcal{A}$ over $M$ can
naturally be regarded as a bundle of Nishimura algebroids over a single point
over $M$.
\end{remark}

In this section an arbitrarily chosen totally intransitive Nishimura algebroid
$\mathcal{A}$\ over $M$ shall be fixed.

\begin{definition}
Given $x\in\mathcal{A}^{\mathcal{D}_{1}}$ and $y\in\mathcal{A}^{\mathcal{D}%
_{2}}$ with $\pi(x)=\pi(y)$, we define $x\circledast y\in\mathcal{A}%
^{\mathcal{D}_{1}\times\mathcal{D}_{2}}$ to be
\[
(d\in\mathcal{D}_{2}\mapsto x)\ast y
\]

\end{definition}

\begin{proposition}
For any $x\in\mathcal{A}^{\mathcal{D}_{1}}$, $y\in\mathcal{A}^{\mathcal{D}%
_{2}}$ and $z\in\mathcal{A}^{\mathcal{D}_{3}}$ with $\pi(x)=\pi(y)=\pi(z) $,
we have
\[
x\circledast(y\circledast z)=(x\circledast y)\circledast z
\]

\end{proposition}

\begin{proof}
This follows simply from the fourth condition in Definition \ref{d2.4}.
\end{proof}

\begin{remark}
By this proposition we can omit parentheses in a combination by $\circledast$.
\end{remark}

The following proposition is the Nishimura algebroid counterpart of
Proposition 3 (\S 3.2) of Lavendhomme \cite{l1}.

\begin{proposition}
\label{tr1}Let $x\in\mathcal{A}^{1}$. Then we have
\begin{align*}
&  \mathcal{A}^{((d_{1},d_{2})\in D(2)\longmapsto d_{1}+d_{2}\in D)}(x)\\
&  =\mathcal{A}^{((d_{1},d_{2})\in D(2)\longmapsto(d_{1},d_{2})\in D^{2}%
)}(x\circledast x)\\
&  =\mathcal{A}^{((d_{1},d_{2})\in D(2)\longmapsto(d_{2},d_{1})\in D^{2}%
)}(x\circledast x)
\end{align*}

\end{proposition}

\begin{proof}
Let $z=\mathcal{A}^{((d_{1},d_{2})\in D(2)\longmapsto(d_{1},d_{2})\in D^{2}%
)}(x\circledast x)$. Then we have
\begin{align*}
&  \mathcal{A}^{(d\in D\longmapsto(d,0)\in D(2))}(z)\\
&  =\mathcal{A}^{(d\in D\longmapsto(d,0)\in D^{2})}(x\circledast x)\\
&  =x
\end{align*}
and
\begin{align*}
&  \mathcal{A}^{(d\in D\longmapsto(0,d)\in D(2))}(z)\\
&  =\mathcal{A}^{(d\in D\longmapsto(0,d)\in D^{2})}(x\circledast x)\\
&  =(d\in D\longmapsto x)(0)\\
&  =x
\end{align*}
Therefore the desired first equality follows at once from the quasi-colimit
diagram in Proposition 6 (\S 2.2) of Lavendhomme \cite{l1}. The desired second
equality can be dealt with similarly.
\end{proof}

The following proposition is the Nishimura algebroid counterpart of
Proposition 6 (\S 3.2) of Lavendhomme \cite{l1}.

\begin{proposition}
\label{tr2}Let $x,y\in\mathcal{A}^{1}$ with $\pi(x)=\pi(y)$. Then we have
\begin{align*}
&  x+y\\
&  =\mathcal{A}^{(d\in D\longmapsto(d,d)\in D^{2})}(y\circledast x)\\
&  =\mathcal{A}^{(d\in D\longmapsto(d,d)\in D^{2})}(x\circledast y)
\end{align*}

\end{proposition}

\begin{proof}
Let $z=\mathcal{A}^{(d_{1},d_{2})\in D(2)\longmapsto(d_{1},d_{2})\in D^{2}%
)}(y\circledast x)$. Then we have
\begin{align*}
&  \mathcal{A}^{(d\in D\longmapsto(d,0)\in D(2))}(z)\\
&  =\mathcal{A}^{(d\in D\longmapsto(d,0)\in D^{2})}(y\circledast x)\\
&  =x
\end{align*}
and
\begin{align*}
&  \mathcal{A}^{(d\in D\longmapsto(0,d)\in D(2))}(z)\\
&  =\mathcal{A}^{(d\in D\longmapsto(0,d)\in D^{2})}(y\circledast x)\\
&  =y
\end{align*}
Therefore it follows from the quasi-colimit diagram in Proposition 6
(\S \ 2.2) of Lavendhomme \cite{l1} that
\begin{align*}
&  x+y\\
&  =\mathcal{A}^{(d\in D\longmapsto(d,d)\in D(2))}(z)\\
&  =\mathcal{A}^{(d\in D\longmapsto(d,d)\in D^{2})}(y\circledast x)
\end{align*}
which establishes the first desired equality. The second desired equality
follows similarly.
\end{proof}

\begin{proposition}
\label{tr3}Given $x,y\in\mathcal{A}^{1}$ with $\pi(x)=\pi(y)$, there exists a
unique $z\in\mathcal{A}^{1}$ with $\pi(x)=\pi(y)=\pi(z)$ such that
\begin{align*}
&  \mathcal{A}^{((d_{1},d_{2})\in D^{2}\longmapsto d_{1}d_{2}\in D)}(z)\\
&  =\mathcal{A}^{((d_{1},d_{2})\in D^{2}\longmapsto(d_{1},d_{2},-d_{1}%
,-d_{2})\in D^{4})}(y\circledast x\circledast y\circledast x)
\end{align*}

\end{proposition}

\begin{proof}
We will show that
\begin{align*}
&  \mathcal{A}^{(d\in D\mapsto(d,0)\in D^{2})}\circ\mathcal{A}^{((d_{1}%
,d_{2})\in D^{2}\longmapsto(d_{1},d_{2},-d_{1},-d_{2})\in D^{4})}(y\circledast
x\circledast y\circledast x)\\
&  =\mathbf{0}_{\pi(x)}^{D}%
\end{align*}
and
\begin{align*}
&  \mathcal{A}^{(d\in D\mapsto(0,d)\in D^{2})}\circ\mathcal{A}^{((d_{1}%
,d_{2})\in D^{2}\longmapsto(d_{1},d_{2},-d_{1},-d_{2})\in D^{4})}(y\circledast
x\circledast y\circledast x)\\
&  =\mathbf{0}_{\pi(y)}^{D}%
\end{align*}
Then the desired result will follow from the quasi-colimit diagram in
Proposition 7 (\S 2.2) of Lavendhomme \cite{l1}. Now we deal with the first
desired identity. Since the composition of $d\in D\mapsto(d,0)\in D^{2}$ and
$(d_{1},d_{2})\in D^{2}\longmapsto(d_{1},d_{2},-d_{1},-d_{2})\in D^{4}$ is
equal to the composition of $d\in D\mapsto(d,d)\in D^{2}$ and $(d_{1}%
,d_{2})\in D^{2}\longmapsto(d_{1},0,-d_{2},0)\in D^{4}$, we have
\begin{align*}
&  \mathcal{A}^{(d\in D\mapsto(d,0)\in D^{2})}\circ\mathcal{A}^{((d_{1}%
,d_{2})\in D^{2}\longmapsto(d_{1},d_{2},-d_{1},-d_{2})\in D^{4})}(y\circledast
x\circledast y\circledast x)\\
&  =\mathcal{A}^{(d\in D\mapsto(d,d)\in D^{2})}\circ\mathcal{A}^{((d_{1}%
,d_{2})\in D^{2}\longmapsto(d_{1},0,-d_{2},0)\in D^{4})}(y\circledast
x\circledast y\circledast x)\\
&  =\mathcal{A}^{(d\in D\mapsto(d,d)\in D^{2})}(\mathcal{A}^{(d_{2}\in
D\longmapsto(-d_{2},0)\in D^{2})}(y\circledast x)\circledast\mathcal{A}%
^{(d_{1}\in D\longmapsto(d_{1},0)\in D^{2})}(y\circledast x))\\
&  =\mathcal{A}^{(d\in D\mapsto(d,d)\in D^{2})}((-x)\circledast x)\\
&  =x-x\text{ \ \ \ \ \ [by Proposition \ref{tr2}]}\\
&  =\mathbf{0}_{\pi(x)}^{D}%
\end{align*}
Now we turn to the second desired identity. Since the composition of $d\in
D\mapsto(0,d)\in D^{2}$ and $(d_{1},d_{2})\in D^{2}\longmapsto(d_{1}%
,d_{2},-d_{1},-d_{2})\in D^{4}$ is equal to the composition of $d\in
D\mapsto(d,d)\in D^{2}$ and $(d_{1},d_{2})\in D^{2}\longmapsto(0,d_{1}%
,0,-d_{2})\in D^{4}$, we have
\begin{align*}
&  \mathcal{A}^{(d\in D\mapsto(0,d)\in D^{2})}\circ\mathcal{A}^{((d_{1}%
,d_{2})\in D^{2}\longmapsto(d_{1},d_{2},-d_{1},-d_{2})\in D^{4})}(y\circledast
x\circledast y\circledast x)\\
&  =\mathcal{A}^{(d\in D\mapsto(d,d)\in D^{2})}\circ\mathcal{A}^{((d_{1}%
,d_{2})\in D^{2}\longmapsto(0,d_{1},0,-d_{2})\in D^{4})}(y\circledast
x\circledast y\circledast x)\\
&  =\mathcal{A}^{(d\in D\mapsto(d,d)\in D^{2})}(\mathcal{A}^{(d_{2}\in
D\longmapsto(0,-d_{2})\in D^{2})}(y\circledast x)\circledast\mathcal{A}%
^{(d_{1}\in D\longmapsto(0,d_{1})\in D^{2})}(y\circledast x))\\
&  =\mathcal{A}^{(d\in D\mapsto(d,d)\in D^{2})}((-y)\circledast y)\\
&  =y-y\text{ \ \ \ \ \ [by Proposition \ref{tr2}]}\\
&  =\mathbf{0}_{\pi(y)}^{D}%
\end{align*}
The proof is now complete.
\end{proof}

\begin{notation}
We will denote the above $z$ by $[x,y]$.
\end{notation}

\begin{proposition}
\label{tr5}Given $x,y\in\mathcal{A}^{1}$ with $\pi(x)=\pi(y)$, we have
\[
\lbrack y,x]=-[x,y]
\]

\end{proposition}

\begin{proof}
Let $m=\pi(x)=\pi(y)$. We have
\begin{align*}
&  \mathcal{A}^{((d_{1},d_{2})\in D^{2}\longmapsto d_{1}d_{2}\in D)}%
\circ\mathcal{A}^{(d\in D\mapsto(d,d)\in D^{2})}([x,y]\circledast\lbrack
y,x])\\
&  =\mathcal{A}^{((d_{1},d_{2})\in D^{2}\longmapsto(d_{1}d_{2},d_{1}d_{2})\in
D^{2})}([x,y]\circledast\lbrack y,x])\\
&  =\mathcal{A}^{((d_{1},d_{2})\in D^{2}\longmapsto(d_{1},d_{2},d_{1}%
,d_{2})\in D^{4})}\circ\mathcal{A}^{((d_{1},d_{2},d_{3},d_{4})\in D^{4}%
\mapsto(d_{1}d_{2},d_{3}d_{4})\in D^{2})}([x,y]\circledast\lbrack y,x])\\
&  =\mathcal{A}^{((d_{1},d_{2})\in D^{2}\longmapsto(d_{1},d_{2},d_{1}%
,d_{2})\in D^{4})}(\mathcal{A}^{((d_{1},d_{2})\in D^{2}\longmapsto d_{1}%
d_{2}\in D)}([x,y])\circledast\mathcal{A}^{((d_{1},d_{2})\in D^{2}\longmapsto
d_{1}d_{2}\in D)}([y,x]))\\
&  =\mathcal{A}^{((d_{1},d_{2})\in D^{2}\longmapsto(d_{1},d_{2},d_{1}%
,d_{2})\in D^{4})}(\mathcal{A}^{((d_{1},d_{2})\in D^{2}\longmapsto d_{1}%
d_{2}\in D)}([x,y])\circledast\\
&  \mathcal{A}^{((d_{1},d_{2})\in D^{2}\longmapsto(d_{2},d_{1})\in D^{2}%
)}\circ\mathcal{A}^{((d_{1},d_{2})\in D^{2}\longmapsto d_{1}d_{2}\in
D)}([y,x]))\\
&  =\mathcal{A}^{((d_{1},d_{2})\in D^{2}\longmapsto(d_{1},d_{2},d_{1}%
,d_{2})\in D^{4})}(\mathcal{A}^{((d_{1},d_{2})\in D^{2}\longmapsto(d_{1}%
,d_{2},-d_{1},-d_{2})\in D^{4})}(y\circledast x\circledast y\circledast
x)\circledast\\
&  \mathcal{A}^{((d_{1},d_{2})\in D^{2}\longmapsto(d_{2},d_{1})\in D^{2}%
)}\circ\mathcal{A}^{((d_{1},d_{2})\in D^{2}\longmapsto(d_{1},d_{2}%
,-d_{1},-d_{2})\in D^{4})}(x\circledast y\circledast x\circledast y))\\
&  =\mathcal{A}^{((d_{1},d_{2})\in D^{2}\longmapsto(d_{1},d_{2},d_{1}%
,d_{2})\in D^{4})}(\mathcal{A}^{((d_{1},d_{2})\in D^{2}\longmapsto(d_{1}%
,d_{2},-d_{1},-d_{2})\in D^{4})}(y\circledast x\circledast y\circledast
x)\circledast\\
&  \mathcal{A}^{((d_{1},d_{2})\in D^{2}\longmapsto(d_{2},d_{1},-d_{2}%
,-d_{1})\in D^{4})}(x\circledast y\circledast x\circledast y))\\
&  =\mathcal{A}^{((d_{1},d_{2})\in D^{2}\longmapsto(d_{1},d_{2},d_{1}%
,d_{2})\in D^{4})}\circ\mathcal{A}^{((d_{1},d_{2},d_{3},d_{4})\in D^{4}%
\mapsto(d_{2},d_{1},-d_{2},-d_{1},d_{3},d_{4},-d_{3},-d_{4})\in D^{8})}\\
&  (y\circledast x\circledast y\circledast x\circledast x\circledast
y\circledast x\circledast y)\\
&  =\mathcal{A}^{((d_{1},d_{2})\in D^{2}\longmapsto(d_{2},d_{1},-d_{2}%
,-d_{1},d_{1},d_{2},-d_{1},-d_{2})\in D^{8})}(y\circledast x\circledast
y\circledast x\circledast x\circledast y\circledast x\circledast y)\\
&  =\mathcal{A}^{((d_{1},d_{2})\in D^{2}\longmapsto(d_{2},d_{1},-d_{2}%
,d_{1},d_{2},-d_{1},-d_{2})\in D^{7})}\circ\mathcal{A}^{((d_{1},d_{2}%
,d_{3},d_{4},d_{5},d_{6},d_{7})\in D^{7}\mapsto(d_{1},d_{2},d_{3},-d_{4}%
,d_{4},d_{5},d_{6},d_{7})\in D^{8})}\\
&  (y\circledast x\circledast y\circledast x\circledast x\circledast
y\circledast x\circledast y)\\
&  =\mathcal{A}^{((d_{1},d_{2})\in D^{2}\longmapsto(d_{2},d_{1},-d_{2}%
,d_{1},d_{2},-d_{1},-d_{2})\in D^{7})}\\
&  (y\circledast x\circledast y\circledast\mathcal{A}^{(d\in D\mapsto(-d,d)\in
D^{2})}(x\circledast x)\circledast y\circledast x\circledast y)\\
&  =\mathcal{A}^{((d_{1},d_{2})\in D^{2}\longmapsto(d_{2},d_{1},-d_{2}%
,d_{2},-d_{1},-d_{2})\in D^{6})}(y\circledast x\circledast y\circledast
y\circledast x\circledast y)\\
&  =\mathcal{A}^{((d_{1},d_{2})\in D^{2}\longmapsto(d_{2},d_{1},d_{2}%
,-d_{1},-d_{2})\in D^{5})}\circ\mathcal{A}^{((d_{1},d_{2},d_{3},d_{4}%
,d_{5})\in D^{5}\longmapsto(d_{1},d_{2},-d_{3},d_{3},d_{4},d_{5})\in D^{6})}\\
&  (y\circledast x\circledast y\circledast y\circledast x\circledast y)\\
&  =\mathcal{A}^{((d_{1},d_{2})\in D^{2}\longmapsto(d_{2},d_{1},d_{2}%
,-d_{1},-d_{2})\in D^{5})}(y\circledast x\circledast\mathcal{A}^{(d\in
D\mapsto(-d,d)\in D^{2})}(y\circledast y)\circledast x\circledast y)\\
&  =\mathcal{A}^{((d_{1},d_{2})\in D^{2}\longmapsto(d_{2},d_{1},-d_{1}%
,-d_{2})\in D^{4})}(y\circledast x\circledast x\circledast y)\\
&  =\mathcal{A}^{((d_{1},d_{2})\in D^{2}\longmapsto(d_{2},d_{1},-d_{2})\in
D^{3})}\circ\mathcal{A}^{((d_{1},d_{2},d_{3})\in D^{3}\longmapsto(d_{1}%
,d_{2},-d_{2},d_{3})\in D^{4})}(y\circledast x\circledast x\circledast y)\\
&  =\mathcal{A}^{((d_{1},d_{2})\in D^{2}\longmapsto(d_{2},d_{1},-d_{2})\in
D^{3})}(y\circledast\mathcal{A}^{(d\in D\mapsto(d,-d)\in D^{2})}(x\circledast
x)\circledast y)\\
&  =\mathcal{A}^{((d_{1},d_{2})\in D^{2}\longmapsto(d_{2},-d_{2})\in D^{2}%
)}(y\circledast y)\\
&  =\mathcal{A}^{((d_{1},d_{2})\in D^{2}\longmapsto d_{2}\in D)}%
\circ\mathcal{A}^{(d\in D\mapsto(d,-d)\in D^{2})}(y\circledast y)\\
&  =\mathbf{0}_{m}^{D^{2}}%
\end{align*}

\end{proof}

\begin{proposition}
\label{tr4}Given $x,y\in\mathcal{A}^{1}$ with $\pi(x)=\pi(y)$, we have
\begin{align}
&  \mathcal{A}^{((d_{1},d_{2})\in D\oplus D\mapsto(d_{1},d_{2})\in D^{2}%
)}(y\circledast x)\nonumber\\
&  =\mathcal{A}^{((d_{1},d_{2})\in D\oplus D\mapsto(d_{1},d_{2})\in D^{2}%
)}\circ\mathcal{A}^{((d_{1},d_{2})\in D^{2}\longmapsto(d_{2},d_{1})\in D^{2}%
)}\mathcal{(}x\circledast y) \label{tr4.3}%
\end{align}
and
\begin{equation}
[x,y]=y\circledast x\overset{\cdot}{-}\mathcal{A}^{((d_{1},d_{2})\in
D^{2}\longmapsto(d_{2},d_{1})\in D^{2})}\mathcal{(}x\circledast y)
\label{tr4.4}%
\end{equation}

\end{proposition}

\begin{proof}
Our proof is the proof of Proposition 8 (\S 3.4) of Lavendhomme \cite{l1} in
disguise. In order to show the identity (\ref{tr4.3}), it suffices, by dint of
the quasi-colimit diagram in Proposition 6 (\S 2.2) of Lavendhomme \cite{l1},
to show that
\begin{align}
&  \mathcal{A}^{(d\in D\longmapsto(d,0)\in D\oplus D)}\circ\mathcal{A}%
^{((d_{1},d_{2})\in D\oplus D\mapsto(d_{1},d_{2})\in D^{2})}(y\circledast
x)\nonumber\\
&  =\mathcal{A}^{(d\in D\longmapsto(d,0)\in D\oplus D)}\circ\mathcal{A}%
^{((d_{1},d_{2})\in D\oplus D\mapsto(d_{1},d_{2})\in D^{2})}\circ\nonumber\\
&  \mathcal{A}^{((d_{1},d_{2})\in D^{2}\longmapsto(d_{2},d_{1})\in D^{2}%
)}\mathcal{(}x\circledast y) \label{tr4.1}%
\end{align}
and
\begin{align}
&  \mathcal{A}^{(d\in D\longmapsto(0,d)\in D\oplus D)}\circ\mathcal{A}%
^{((d_{1},d_{2})\in D\oplus D\mapsto(d_{1},d_{2})\in D^{2})}(y\circledast
x)\nonumber\\
&  =\mathcal{A}^{(d\in D\longmapsto(0,d)\in D\oplus D)}\circ\mathcal{A}%
^{((d_{1},d_{2})\in D\oplus D\mapsto(d_{1},d_{2})\in D^{2})}\circ\nonumber\\
&  \mathcal{A}^{((d_{1},d_{2})\in D^{2}\longmapsto(d_{2},d_{1})\in D^{2}%
)}\mathcal{(}x\circledast y) \label{tr4.2}%
\end{align}
Since the composition of $d\in D\longmapsto(d,0)\in D\oplus D$ and
$(d_{1},d_{2})\in D\oplus D\mapsto(d_{1},d_{2})\in D^{2}$ is equal to $d\in
D\longmapsto(d,0)\in D^{2}$, and since the composition of $d\in D\longmapsto
(d,0)\in D\oplus D$, $(d_{1},d_{2})\in D\oplus D\mapsto(d_{1},d_{2})\in D^{2}$
and $(d_{1},d_{2})\in D^{2}\longmapsto(d_{2},d_{1})\in D^{2}$ is equal to
$d\in D\longmapsto(0,d)\in D^{2}$, it is easy to see that both sides of the
identity (\ref{tr4.1}) are equal to $x$ by the second condition in Definition
\ref{d2.4}. The identity (\ref{tr4.2}) can be established similarly. Let
\[
z=\mathcal{A}^{((d_{1},d_{2},d_{3})\in D^{2}\oplus D\longmapsto(d_{2}%
,d_{3},d_{1})\in D^{3})}\mathcal{(}x\circledast\lbrack x,y]\circledast y)
\]
Then we have
\begin{align*}
&  \mathcal{A}^{((d_{1},d_{2})\in D^{2}\longmapsto(d_{1},d_{2},0)\in
D^{2}\oplus D)}(z)\\
&  =\mathcal{A}^{((d_{1},d_{2})\in D^{2}\longmapsto(d_{1},d_{2},0)\in
D^{2}\oplus D)}\circ\mathcal{A}^{((d_{1},d_{2},d_{3})\in D^{2}\oplus
D\longmapsto(d_{2},d_{3},d_{1})\in D^{3})}\mathcal{(}x\circledast\lbrack
x,y]\circledast y)\\
&  =\mathcal{A}^{((d_{1},d_{2})\in D^{2}\longmapsto(d_{2},0,d_{1})\in D^{3}%
)}\mathcal{(}x\circledast\lbrack x,y]\circledast y)\\
&  =\mathcal{A}^{((d_{1},d_{2})\in D^{2}\longmapsto(d_{2},d_{1})\in D^{2}%
)}\circ\mathcal{A}^{((d_{1},d_{2})\in D^{2}\longmapsto(d_{1},0,d_{2})\in
D^{3})}\mathcal{(}x\circledast\lbrack x,y]\circledast y)\\
&  =\mathcal{A}^{((d_{1},d_{2})\in D^{2}\longmapsto(d_{2},d_{1})\in D^{2}%
)}(x\circledast\mathcal{A}^{(d\in D\longmapsto(d,0)\in D^{2})}%
([x,y]\circledast y))\\
&  =\mathcal{A}^{((d_{1},d_{2})\in D^{2}\longmapsto(d_{2},d_{1})\in D^{2}%
)}(x\circledast y)
\end{align*}
while we have
\begin{align*}
&  \mathcal{A}^{((d_{1},d_{2})\in D^{2}\longmapsto(d_{1},d_{2},d_{1}d_{2})\in
D^{2}\oplus D)}(z)\\
&  =\mathcal{A}^{((d_{1},d_{2})\in D^{2}\longmapsto(d_{1},d_{2},d_{1}d_{2})\in
D^{2}\oplus D)}\circ\mathcal{A}^{((d_{1},d_{2},d_{3})\in D^{2}\oplus
D\longmapsto(d_{2},d_{3},d_{1})\in D^{3})}\mathcal{(}x\circledast\lbrack
x,y]\circledast y)\\
&  =\mathcal{A}^{((d_{1},d_{2})\in D^{2}\longmapsto(d_{2},d_{1}d_{2},d_{1})\in
D^{3})}\mathcal{(}x\circledast\lbrack x,y]\circledast y)\\
&  =\mathcal{A}^{((d_{1},d_{2})\in D^{2}\longmapsto(d_{2},d_{1}d_{2},d_{1})\in
D^{3})}\mathcal{(}x\circledast\mathcal{A}^{(d\in D\mapsto-d\in D)}%
\circ\mathcal{A}^{(d\in D\mapsto-d\in D)}([x,y])\circledast y)\\
&  =\mathcal{A}^{((d_{1},d_{2})\in D^{2}\longmapsto(d_{2},d_{1}d_{2},d_{1})\in
D^{3})}\circ\mathcal{A}^{((d_{1},d_{2},d_{3})\in D^{3}\mapsto(d_{1}%
,-d_{2},d_{3})\in D^{3})}(x\circledast\mathcal{A}^{(d\in D\mapsto-d\in
D)}([x,y])\circledast y)\\
&  =\mathcal{A}^{((d_{1},d_{2})\in D^{2}\longmapsto(d_{2},-d_{1}d_{2}%
,d_{1})\in D^{3})}(x\circledast\lbrack y,x]\circledast y)\\
&  \text{[By Proposition \ref{tr5}]}\\
&  =\mathcal{A}^{((d_{1},d_{2})\in D^{2}\longmapsto(d_{2},-d_{2},d_{1}%
,d_{1})\in D^{4})}\circ\mathcal{A}^{((d_{1},d_{2},d_{3},d_{4})\in D^{4}%
\mapsto(d_{1},d_{2}d_{3},d_{4})\in D^{3})}(x\circledast\lbrack y,x]\circledast
y)\\
&  =\mathcal{A}^{((d_{1},d_{2})\in D^{2}\longmapsto(d_{2},-d_{2},d_{1}%
,d_{1})\in D^{4})}(x\circledast\mathcal{A}^{((d_{1},d_{2})\in D^{2}\mapsto
d_{1}d_{2}\in D)}([y,x])\circledast y)\\
&  =\mathcal{A}^{((d_{1},d_{2})\in D^{2}\longmapsto(d_{2},-d_{2},d_{1}%
,d_{1})\in D^{4})}(x\circledast\mathcal{A}^{((d_{1},d_{2})\in D^{2}%
\longmapsto(d_{1},d_{2},-d_{1},-d_{2})\in D^{4})}(x\circledast y\circledast
x\circledast y)\circledast y)\\
&  =\mathcal{A}^{((d_{1},d_{2})\in D^{2}\longmapsto(d_{2},-d_{2},d_{1}%
,d_{1})\in D^{4})}\circ\mathcal{A}^{((d_{1},d_{2},d_{3},d_{4})\in D^{4}%
\mapsto(d_{1},d_{2},d_{3},-d_{2},-d_{3},d_{4})\in D^{6}})\\
&  (x\circledast x\circledast y\circledast x\circledast y\circledast y)\\
&  =\mathcal{A}^{((d_{1},d_{2})\in D^{2}\longmapsto(d_{2},-d_{2},d_{1}%
,d_{2},-d_{1},d_{1})\in D^{6})}(x\circledast x\circledast y\circledast
x\circledast y\circledast y)\\
&  =\mathcal{A}^{((d_{1},d_{2})\in D^{2}\longmapsto(d_{2},d_{1},d_{2}%
,d_{1})\in D^{4})}\circ\mathcal{A}^{((d_{1},d_{2},d_{3},d_{4})\in
D^{4}\longmapsto(d_{1},-d_{1},d_{2},d_{3},-d_{4},d_{4})\in D^{6})}\\
&  (x\circledast x\circledast y\circledast x\circledast y\circledast y)\\
&  =\mathcal{A}^{((d_{1},d_{2})\in D^{2}\longmapsto(d_{2},d_{1},d_{2}%
,d_{1})\in D^{4})}(\mathcal{A}^{(d\in D\mapsto(-d,d)\in D^{2})}\mathcal{(}%
x\circledast x)\circledast y\circledast x\circledast\\
&  \mathcal{A}^{(d\in D\mapsto(d,-d)\in D^{2})}\mathcal{(}y\circledast y))\\
&  =\mathcal{A}^{((d_{1},d_{2})\in D^{2}\longmapsto(d_{2},d_{1},d_{2}%
,d_{1})\in D^{4})}(\mathbf{0}_{m}^{D}\circledast y\circledast x\circledast
\mathbf{0}_{m}^{D})\\
&  =y\circledast x
\end{align*}
Therefore we have
\begin{align*}
&  y\circledast x\overset{\cdot}{-}\mathcal{A}^{((d_{1},d_{2})\in
D^{2}\longmapsto(d_{2},d_{1})\in D^{2})}\mathcal{(}x\circledast y)\\
&  =\mathcal{A}^{(d\in D\mapsto(0,0,d)\in D^{2}\oplus D)}(z)\\
&  =\mathcal{A}^{(d\in D\mapsto(0,0,d)\in D^{2}\oplus D)}\circ\mathcal{A}%
^{((d_{1},d_{2},d_{3})\in D^{2}\oplus D\longmapsto(d_{2},d_{3},d_{1})\in
D^{3})}\mathcal{(}x\circledast\lbrack x,y]\circledast y)\\
&  =[x,y]
\end{align*}
This completes the proof.
\end{proof}

\begin{proposition}
\begin{enumerate}
\item Given $x\in\mathcal{A}^{1}$ and $y,z\in\mathcal{A}^{2}$with $\pi
(x)=\pi(y)=\pi(z)$, if we have
\[
\mathcal{A}^{((d_{1},d_{2})\in D\oplus D\mapsto(d_{1},d_{2})\in D^{2}%
)}(y)=\mathcal{A}^{((d_{1},d_{2})\in D\oplus D\mapsto(d_{1},d_{2})\in D^{2}%
)}(z)
\]
then we have
\[
\mathcal{A}^{((d_{1},d_{2},d_{3})\in D\times(D\oplus D)\mapsto(d_{1}%
,d_{2},d_{3})\in D^{3})}(y\circledast x)=\mathcal{A}^{((d_{1},d_{2},d_{3})\in
D\times(D\oplus D)\mapsto(d_{1},d_{2},d_{3})\in D^{3})}(z\circledast x)
\]
and
\[
z\circledast x\overset{\cdot}{\underset{1}{-}}y\circledast x=(z\overset{\cdot
}{-}y)\circledast x
\]

\item Given $x,y\in\mathcal{A}^{2}$ and $z\in\mathcal{A}^{1}$with $\pi
(x)=\pi(y)=\pi(z)$, if we have
\[
\mathcal{A}^{((d_{1},d_{2})\in D\oplus D\mapsto(d_{1},d_{2})\in D^{2}%
)}(x)=\mathcal{A}^{((d_{1},d_{2})\in D\oplus D\mapsto(d_{1},d_{2})\in D^{2}%
)}(y)
\]
then we have
\[
\mathcal{A}^{((d_{1},d_{2},d_{3})\in(D\oplus D)\times D\mapsto(d_{1}%
,d_{2},d_{3})\in D^{3})}(z\circledast x)=\mathcal{A}^{((d_{1},d_{2},d_{3}%
)\in(D\oplus D)\times D\mapsto(d_{1},d_{2},d_{3})\in D^{3})}(z\circledast y)
\]
and
\[
\mathcal{A}^{((d_{1},d_{2})\in D^{2}\mapsto(d_{2},d_{1})\in D^{2}%
)}(z\circledast y\overset{\cdot}{\underset{3}{-}}z\circledast x)=z\circledast
(y\overset{\cdot}{-}x)
\]

\end{enumerate}
\end{proposition}

\begin{proof}
This follows simply from Proposition \ref{t2.6}.
\end{proof}

\begin{proposition}
Given $x,y,z\in\mathcal{A}^{1}$ with $\pi(x)=\pi(y)=\pi(z)$, let it be the
case that
\begin{align*}
u_{123}  &  =z\circledast y\circledast x\\
u_{132}  &  =\mathcal{A}^{((d_{1},d_{2},d_{3})\in D^{3}\mapsto(d_{1}%
,d_{3},d_{2})\in D^{3})}\mathcal{(}y\circledast z\circledast x)\\
u_{213}  &  =\mathcal{A}^{((d_{1},d_{2},d_{3})\in D^{3}\mapsto(d_{2}%
,d_{1},d_{3})\in D^{3})}(z\circledast x\circledast y)\\
u_{231}  &  =\mathcal{A}^{((d_{1},d_{2},d_{3})\in D^{3}\mapsto(d_{2}%
,d_{3},d_{1})\in D^{3})}(x\circledast z\circledast y)\\
u_{312}  &  =\mathcal{A}^{((d_{1},d_{2},d_{3})\in D^{3}\mapsto(d_{3}%
,d_{1},d_{2})\in D^{3})}\mathcal{(}y\circledast x\circledast z)\\
u_{321}  &  =\mathcal{A}^{((d_{1},d_{2},d_{3})\in D^{3}\mapsto(d_{3}%
,d_{2},d_{1})\in D^{3})}\mathcal{(}x\circledast y\circledast z)
\end{align*}
Then the right-hands of the following three identities are meaningful, and all
the three identities hold:
\begin{align*}
\lbrack x,[y,z]]  &  =(u_{123}\overset{\cdot}{\underset{1}{-}}u_{132}%
)\overset{\cdot}{-}(u_{231}\overset{\cdot}{\underset{1}{-}}u_{321})\\
\lbrack y,[z,x]]  &  =(u_{231}\overset{\cdot}{\underset{2}{-}}u_{213}%
)\overset{\cdot}{-}(u_{312}\overset{\cdot}{\underset{2}{-}}u_{132})\\
\lbrack z,[x,y]]  &  =(u_{312}\overset{\cdot}{\underset{3}{-}}u_{321}%
)\overset{\cdot}{-}(u_{123}\overset{\cdot}{\underset{3}{-}}u_{213})
\end{align*}

\end{proposition}

\begin{proof}
Here we deal only with the first identity, leaving the other two identities to
the reader. We have
\begin{align*}
&  [x,[y,z]]\\
&  =[y,z]\circledast x\overset{\cdot}{-}\mathcal{A}^{((d_{1},d_{2})\in
D^{2}\longmapsto(d_{2},d_{1})\in D^{2})}\mathcal{(}x\circledast\lbrack y,z])\\
&  =\{z\circledast y\overset{\cdot}{-}\mathcal{A}^{((d_{1},d_{2})\in
D^{2}\longmapsto(d_{2},d_{1})\in D^{2})}\mathcal{(}y\circledast
z)\}\circledast x\\
&  \overset{\cdot}{-}\mathcal{A}^{((d_{1},d_{2})\in D^{2}\longmapsto
(d_{2},d_{1})\in D^{2})}\mathcal{(}x\circledast\{z\circledast y\overset{\cdot
}{-}\mathcal{A}^{((d_{1},d_{2})\in D^{2}\longmapsto(d_{2},d_{1})\in D^{2}%
)}\mathcal{(}y\circledast z)\})\\
&  =\{z\circledast y\circledast x\overset{\cdot}{\underset{1}{-}}%
\mathcal{A}^{((d_{1},d_{2})\in D^{2}\longmapsto(d_{2},d_{1})\in D^{2}%
)}\mathcal{(}y\circledast z)\circledast x\}\\
&  \overset{\cdot}{-}\mathcal{\{}x\circledast z\circledast y\overset{\cdot
}{\underset{3}{-}}x\circledast\mathcal{A}^{((d_{1},d_{2})\in D^{2}%
\longmapsto(d_{2},d_{1})\in D^{2})}\mathcal{(}y\circledast z)\}\\
&  =\{z\circledast y\circledast x\overset{\cdot}{\underset{1}{-}}%
\mathcal{A}^{((d_{1},d_{2},d_{3})\in D^{3}\mapsto(d_{1},d_{3},d_{2})\in
D^{3})}\mathcal{(}y\circledast z\circledast x)\}\\
&  \overset{\cdot}{-}\mathcal{\{A}^{((d_{1},d_{2},d_{3})\in D^{3}\mapsto
(d_{2},d_{3},d_{1})\in D^{3})}(x\circledast z\circledast y)\\
&  \overset{\cdot}{\underset{1}{-}}\mathcal{A}^{((d_{1},d_{2},d_{3})\in
D^{3}\mapsto(d_{2},d_{3},d_{1})\in D^{3})}(x\circledast\mathcal{A}%
^{((d_{1},d_{2})\in D^{2}\longmapsto(d_{2},d_{1})\in D^{2})}\mathcal{(}%
y\circledast z))\}\\
&  =\{z\circledast y\circledast x\overset{\cdot}{\underset{1}{-}}%
\mathcal{A}^{((d_{1},d_{2},d_{3})\in D^{3}\mapsto(d_{1},d_{3},d_{2})\in
D^{3})}\mathcal{(}y\circledast z\circledast x)\}\\
&  \overset{\cdot}{-}\mathcal{\{A}^{((d_{1},d_{2},d_{3})\in D^{3}\mapsto
(d_{2},d_{3},d_{1})\in D^{3})}(x\circledast z\circledast y)\\
&  \overset{\cdot}{\underset{1}{-}}\mathcal{A}^{((d_{1},d_{2},d_{3})\in
D^{3}\mapsto(d_{2},d_{3},d_{1})\in D^{3})}\circ\mathcal{A}^{((d_{1}%
,d_{2},d_{3})\in D^{3}\longmapsto(d_{2},d_{1},d_{3})\in D^{3})}(x\circledast
y\circledast z)\}\\
&  =\{z\circledast y\circledast x\overset{\cdot}{\underset{1}{-}}%
\mathcal{A}^{((d_{1},d_{2},d_{3})\in D^{3}\mapsto(d_{1},d_{3},d_{2})\in
D^{3})}\mathcal{(}y\circledast z\circledast x)\}\\
&  \overset{\cdot}{-}\mathcal{\{A}^{((d_{1},d_{2},d_{3})\in D^{3}\mapsto
(d_{2},d_{3},d_{1})\in D^{3})}(x\circledast z\circledast y)\overset{\cdot
}{\underset{1}{-}}\mathcal{A}^{((d_{1},d_{2},d_{3})\in D^{3}\mapsto
(d_{3},d_{2},d_{1})\in D^{3})}\mathcal{(}x\circledast y\circledast z)\}\\
&  =(u_{123}\overset{\cdot}{\underset{1}{-}}u_{132})\overset{\cdot}{-}%
(u_{231}\overset{\cdot}{\underset{1}{-}}u_{321})
\end{align*}

\end{proof}

\begin{theorem}
Given $m\in M$, the Jacobi identity holds for $\mathcal{A}_{m}^{1}$ with
respect to the Lie bracket $[\cdot,\cdot]$. I.e., we have
\[
\lbrack x,[y,z]]+[y,[z,x]]+[z,[x,y]]=\mathbf{0}%
\]
for any $x,y,z\in\mathcal{A}_{m}^{1}$.
\end{theorem}

\section{From Nishimura Algebroids to Lie Algebroids}

Let $\mathcal{A}$ be a \textit{Nishimura algebroid over }$M$. It is very easy
to see that

\begin{proposition}
\label{lt2}By assigning $\Gamma(\mathcal{A)}^{\mathcal{D}}\mathcal{=}%
\Gamma(\mathcal{A}^{\mathcal{D}}\mathcal{)}$ to each object $\mathcal{D}$ in
$\mathbf{Simp}$ and assigning $\Gamma(\mathcal{A)}^{f}:\gamma\in$
$\Gamma(\mathcal{A}^{\mathcal{D}^{\prime}}\mathcal{)\mapsto A}^{f}\circ
\gamma\in\Gamma(\mathcal{A}^{\mathcal{D}}\mathcal{)}$ to each morphism
$f:\mathcal{D\rightarrow D}^{\prime}$ in $\mathbf{Simp}$, we have a
\textit{Nishimura algebroid}$_{2}$ $\Gamma(\mathcal{A)}$ \textit{over} a
single point, where $\Gamma(\mathcal{A}^{\mathcal{D}}\mathcal{)}$ denotes the
space of global sections of the bundle $\mathcal{A}^{\mathcal{D}}$ over $M$.
Endowed with the trivial anchor natural transformation, it is a
\textit{Nishimura algebroid}$_{3}$ \textit{over} a single point.
\end{proposition}

\begin{definition}
Given $X\in\Gamma(\mathcal{A}^{\mathcal{D}_{1}}\mathcal{)}$ and $Y\in
\Gamma(\mathcal{A}^{\mathcal{D}_{2}}\mathcal{)}$, we define $Y\circledcirc
X\in\Gamma(\mathcal{A}^{\mathcal{D}_{1}\times\mathcal{D}_{2}}\mathcal{)}$ to
be
\[
(Y\circledcirc X)_{m}=(Y\circ\mathbf{a}(X_{m}))\ast X_{m}%
\]
for any $m\in M$.
\end{definition}

Now we have

\begin{proposition}
\label{lt1}Given $X\in\Gamma(\mathcal{A}^{\mathcal{D}_{1}}\mathcal{)}$,
$Y\in\Gamma(\mathcal{A}^{\mathcal{D}_{2}}\mathcal{)}$ and $Z\in\Gamma
(\mathcal{A}^{\mathcal{D}_{3}}\mathcal{)}$, we have
\[
Z\circledcirc(Y\circledcirc X)=(Z\circledcirc Y)\circledcirc X
\]

\end{proposition}

\begin{proof}
Let $m\in M$. We have
\begin{align*}
&  (Z\circledcirc(Y\circledcirc X))_{m}\\
&  =(Z\circ\mathbf{a}((Y\circledcirc X)_{m}))\ast(Y\circledcirc X)_{m}\\
&  =(Z\circ(\mathbf{a}(Y)\circledcirc\mathbf{a}(X))_{m})\ast\{(Y\circ
\mathbf{a}(X_{m}))\ast X_{m}\}\\
&  =[\{m^{\prime}\in M\mapsto(Z\circ\mathbf{a}(Y_{m^{\prime}}))\ast
Y_{m^{\prime}}\}\circ\mathbf{a}(X_{m})]\ast X_{m}\\
&  \text{[By the fourth condition in Definition \ref{d2.4}]}\\
&  =((Z\circledcirc Y)\circledcirc X)_{m}%
\end{align*}

\end{proof}

\begin{remark}
By this proposition we can omit parentheses in a combination by $\circledcirc$.
\end{remark}

\begin{proposition}
\label{lt3}By adopting $\circledcirc$ as $\ast_{\Gamma(\mathcal{A)}}$, our
\textit{Nishimura algebroid}$_{3}$ $\Gamma(\mathcal{A)}$ \textit{over} a
single point is a \textit{Nishimura algebroid}$_{4}$ \textit{over} a single point.
\end{proposition}

\begin{proof}
The fourth condition in Definition \ref{d2.4} follows from Proposition
\ref{lt1}. The other three conditions follow trivially.
\end{proof}

Therefore all the discussions of the previous section hold. In particular, we have

\begin{theorem}
\label{lt4}Given $X,Y\in\Gamma(\mathcal{A}^{1}\mathcal{)}$, we can define
$[X,Y]\in\Gamma(\mathcal{A}^{1}\mathcal{)}$ to be the unique one satisfying
\begin{align*}
&  \mathcal{A}^{((d_{1},d_{2})\in D^{2}\longmapsto d_{1}d_{2}\in D)}%
\circ\lbrack X,Y]\\
&  =\mathcal{A}^{((d_{1},d_{2})\in D^{2}\longmapsto(d_{1},d_{2},-d_{1}%
,-d_{2})\in D^{4})}\circ(Y\circledcirc X\circledcirc Y\circledcirc X)\text{,}%
\end{align*}
with respect to which $\Gamma(\mathcal{A}^{1}\mathcal{)}$ is a Lie algebra.
\end{theorem}

\begin{proposition}
\label{lt5}Given $X,Y\in\Gamma(\mathcal{A}^{1}\mathcal{)}$ and $f\in
\mathbb{R}^{M}$, we have
\[
Y\circledcirc fX=f\underset{1}{\cdot}(Y\circledcirc X)
\]

\end{proposition}

\begin{proof}
Let $m\in M$.
\begin{align*}
&  (Y\circledcirc fX)_{m}\\
&  =(Y\circ\mathbf{a}(f(m)X_{m}))\ast(f(m)X_{m})\\
&  =(Y\circ(f(m)\mathbf{a}(X_{m})))\ast(f(m)X_{m})\\
&  =f(m)\underset{1}{\cdot}(Y\circledcirc X)_{m}%
\end{align*}

\end{proof}

\begin{proposition}
\label{lt6}Given $X,Y\in\Gamma(\mathcal{A}^{1}\mathcal{)}$ and $f\in
\mathbb{R}^{M}$, we have
\[
fY\circledcirc X\overset{\cdot}{-}f\underset{2}{\cdot}(Y\circledcirc
X)=\mathbf{a}(X)(f)Y
\]

\end{proposition}

\begin{proof}
Let $m\in M$. We define $\mu\in\mathcal{A}_{m}^{D^{2}\oplus D}$ to be
\[
\mu=\mathcal{A}^{((d_{1},d_{2},d_{3})\in D^{2}\oplus D\mapsto(d_{1}%
,d_{2}f(m)+d_{3}\mathbf{a}(X_{m})(f))\in D^{2})}((Y\circledcirc X)_{m})
\]
where $\mathbf{a}(X_{m})(f)$ is the Lie derivative of $f$ with respect
$\mathbf{a}(X_{m})$. It is easy to see that
\begin{align*}
&  \mathcal{A}^{((d_{1},d_{2})\in D^{2}\mapsto(d_{1},d_{2},d_{1}d_{2})\in
D^{2}\oplus D)}(\mu)\\
&  =\mathcal{A}^{((d_{1},d_{2})\in D^{2}\mapsto(d_{1},d_{2},d_{1}d_{2})\in
D^{2}\oplus D)}\circ\mathcal{A}^{((d_{1},d_{2},d_{3})\in D^{2}\oplus
D\mapsto(d_{1},d_{2}f(m)+d_{3}\mathbf{a}(X_{m})(f))\in D^{2})}((Y\circledcirc
X)_{m})\\
&  =\mathcal{A}^{((d_{1},d_{2})\in D^{2}\mapsto(d_{1},d_{2}f(m)+d_{1}%
d_{2}\mathbf{a}(X_{m})(f))\in D^{2})}((Y\circledcirc X)_{m})\\
&  =\mathcal{A}^{((d_{1},d_{2})\in D^{2}\mapsto(d_{1},d_{2}f(m)+d_{1}%
d_{2}\mathbf{a}(X_{m})(f))\in D^{2})}((Y\circledcirc X)_{m})\\
&  =\mathcal{A}^{((d_{1},d_{2})\in D^{2}\mapsto(d_{1},d_{2}f(\mathbf{a}%
(X_{m})(d_{1})))\in D^{2})}((Y\circledcirc X)_{m})\\
&  =(fY\circledcirc X)_{m}\\
&  \text{\lbrack By the third condition in Definition \ref{d2.4}]}%
\end{align*}
It is also easy to see that
\begin{align*}
&  \mathcal{A}^{((d_{1},d_{2})\in D^{2}\mapsto(d_{1},d_{2},0)\in D^{2}\oplus
D)}(\mu)\\
&  =\mathcal{A}^{((d_{1},d_{2})\in D^{2}\mapsto(d_{1},d_{2},0)\in D^{2}\oplus
D)}\circ\mathcal{A}^{((d_{1},d_{2},d_{3})\in D^{2}\oplus D\mapsto(d_{1}%
,d_{2}f(m)+d_{3}\mathbf{a}(X_{m})(f))\in D^{2})}((Y\circledcirc X)_{m})\\
&  =\mathcal{A}^{((d_{1},d_{2})\in D^{2}\mapsto(d_{1},d_{2}f(m))\in D^{2}%
)}((Y\circledcirc X)_{m})\\
&  =(f\underset{2}{\cdot}(Y\circledcirc X))_{m}\\
&  \text{\lbrack By the third condition in Definition \ref{d2.4}]}%
\end{align*}
Therefore we have
\begin{align*}
&  (fY\circledcirc X)_{m}\overset{\cdot}{-}(f\underset{2}{\cdot}(Y\circledcirc
X))_{m}\\
&  =\mathcal{A}^{(d\in D\mapsto(0,0,d)\in D^{2}\oplus D)}(\mu)\\
&  =\mathcal{A}^{(d\in D\mapsto(0,0,d)\in D^{2}\oplus D)}\circ\mathcal{A}%
^{((d_{1},d_{2},d_{3})\in D^{2}\oplus D\mapsto(d_{1},d_{2}f(m)+d_{3}%
\mathbf{a}(X_{m})(f))\in D^{2})}((Y\circledcirc X)_{m})\\
&  =\mathcal{A}^{(d\in D\mapsto(0,d\mathbf{a}(X_{m})(f))\in D^{2}%
)}((Y\circledcirc X)_{m})\\
&  =\mathbf{a}(X_{m})(f)Y_{m}%
\end{align*}
This completes the proof.
\end{proof}

\begin{proposition}
\label{lt7}Given $X,Y\in\Gamma(\mathcal{A}^{1}\mathcal{)}$ and $f\in
\mathbb{R}^{M}$, we have
\[
\lbrack X,fY]=f[X,Y]+\mathbf{a}(X)(f)Y
\]

\end{proposition}

\begin{proof}
We have
\begin{align*}
&  [X,fY]\\
&  =fY\circledcirc X\overset{\cdot}{-}\mathcal{A}^{((d_{1},d_{2})\in
D^{2}\longmapsto(d_{2},d_{1})\in D^{2})}\mathcal{(}X\circledcirc fY)\\
&  =\{fY\circledcirc X\overset{\cdot}{-}f\underset{2}{\cdot}(Y\circledcirc
X)\}-\\
&  \{\mathcal{A}^{((d_{1},d_{2})\in D^{2}\longmapsto(d_{2},d_{1})\in D^{2}%
)}(f\underset{1}{\cdot}(X\circledcirc Y))\overset{\cdot}{-}f\underset{2}%
{\cdot}(Y\circledcirc X)\}\\
&  \text{[By Proposition \ref{lt6}]}\\
&  =\mathbf{a}(X)(f)Y+f[X,Y]
\end{align*}
This completes the proof.
\end{proof}

\begin{theorem}
Given a \textit{Nishimura algebroid }$\mathcal{A}$ \textit{over }$M$,
$\mathcal{A}^{1}$ is a Lie algebroid over $M$.
\end{theorem}

\begin{proof}
This follows directly from Theorem \ref{lt4} and Proposition \ref{lt7}.
\end{proof}

\end{document}